\documentclass[10pt]{amsart}
\usepackage{latexsym}
\usepackage{amsmath}
\usepackage{amssymb}
\usepackage{epsfig}
\usepackage[matrix,arrow,cmtip,curve,dvips]{xy}

\begin{document}

\theoremstyle{plain}
    \newtheorem{thm}{Theorem}[section]
    \newtheorem{prop}[thm]{Proposition}
    \newtheorem{lemma}[thm]{Lemma}
    \newtheorem{conj}[thm]{Conjecture}
    \newtheorem{cor}[thm]{Corollary}

\theoremstyle{definition}
    \newtheorem{defn}[thm]{Definition}

\theoremstyle{remark}
    \newtheorem{rem}[thm]{Remark}
    \newtheorem{example}[thm]{Example}
    \newtheorem{question}[thm]{Question}

\newcommand{\Poset}{\mbox{\upshape Poset}}
\newcommand{\ind}{\mbox{\upshape ind}}
\newcommand{\Hom}{\mbox{\upshape Hom}}
\newcommand{\conn}{\mbox{\upshape conn}}
\newcommand{\Ht}{\mbox{\upshape ht}}
\newcommand{\colim}{\mbox{\upshape colim}}
\newcommand{\Max}{\mbox{\upshape max}}
\newcommand{\Min}{\mbox{\upshape min}}
\newcommand{\Aut}{\mbox{\upshape Aut}}
\newcommand{\diam}{\mbox{\upshape diam}}
\newcommand{\bd}{\mbox{\upshape bd}}

% Article information
\title{The universality of Hom complexes}
\date{\today}

% Author information
\author{Anton Dochtermann}
\address {Department of Mathematics \\
Box 354350 \\
          University of Washington \\
          Seattle, WA. 98195 } \email{antondoc@math.washington.edu}
\urladdr{http://www.math.washington.edu/$\sim$antondoc}

% AMS information
% \keywords{Graphs, graph homomorphism, Hom complex}
% \subjclass[2000]{Primary: 05C75, 57M15; Secondary: 55U35, 18G35}

\begin{abstract}
It is shown that if $T$ is a connected nontrivial graph and $X$ is an arbitrary finite simplicial complex, then there is a graph $G$ such
that the complex $\Hom(T,G)$ is homotopy equivalent to $X$.  The
proof is constructive, and uses a nerve lemma.  Along the way
several results regarding $\Hom$ complexes, exponentials, and
subdivision are established that may be of independent interest.
\end{abstract}

\maketitle

\section{Introduction}
The $\Hom$ complex is a functorial way to assign a poset (and hence
topological space) $\Hom(T,G)$ to a pair of graphs $T$ and $G$.
Versions of these spaces were introduced by Lov\'{a}sz in his proof
of Kneser's conjecture (~\cite{Lov78}), and later further
investigated by Babson and Kozlov in ~\cite{BKcom} and
~\cite{BKpro}. The automorphism group of $T$ naturally acts on the
space $\Hom(T,G)$, and in the case that $T = K_2$ is an edge and $G$
is graph without loops, the complex $\Hom(T,G)$ is a space with a
\textit{free} ${\mathbb Z}_2$-action. In ~\cite{Cs05} Csorba shows
that \textit{any} free ${\mathbb Z}_2$-space can be realized (up to
${\mathbb Z}_2$-homotopy type) as $\Hom(K_2, G)$ for some suitably
chosen graph $G$.  His proof involves a simple and elegant
construction in which one obtains a graph $G$ whose vertices are
precisely those of the given ${\mathbb Z}_2$-simplicial complex.

A natural question to ask is what homotopy types can be realized as
$\Hom(T,?)$ for other test graphs $T$. As Csorba points out,
arbitrary homotopy types cannot be realized by $\Hom$ complexes of
\textit{loopless} graphs even with $T=K_2$ as the test graph; all
such $\Hom$ complexes will be free ${\mathbb Z}_2$-spaces and hence
will present topological obstructions (e.g. parity of the Euler
characteristic). However, if we allow loops on our graphs, and do
not concern ourselves with group actions, we are able to prove the
following `universality' of $\Hom$ complexes.

\begin{thm} \label{maintheorem}
Let $T$ be a connected graph with at least one edge, and suppose $X$
is a finite simplicial complex.  Then there exists a graph $G_{k,X}$
(depending on $X$ and the diameter of $T$) and a homotopy
equivalence

\begin{center}
$\Hom(T,G_{k,X}) \simeq X$.
\end{center}
\end{thm}

The graph $G_{k,X}$ will be \textit{reflexive} (that is, has loops
on all the vertices), and hence the space $\Hom(T,G_{k,X})$ will no
longer carry a free $\Aut(T)$ action.  The idea behind our proof of
this theorem will be to consider $X^k = \bd^k(X)$, a high enough
(depending on the diameter of $T$) barycentric subdivision of the
given simplicial complex $X$, and to define $G_{k,X}$ as the
1-skeleton of $X^k$ with loops placed on each vertex. To show that
$\Hom(T,G_{k,X})$ has the desired homotopy type, we will first
replace it with a homotopy equivalent space $X^\prime$ (which will
be the clique complex of some graph). We then determine the homotopy
type of $X^\prime$ by covering it with a collection of contractible
subcomplexes (with contractible intersections) and then employing a
nerve lemma.

The structure of the paper is as follows.  In section 2 we provide
some necessary background on graphs, $\Hom$ complexes, and their
properties.  Section 3 is devoted to the proof of the main result
and some related lemmas.  We conclude in section 4 with some open
questions.

\section{Main objects of study}

In this section we record some basic facts about graphs and $\Hom$
complexes.  For us, a \textit{graph} $G = (V(G), E(G))$ consists of
a vertex set $V(G)$ and an edge set $E(G) \subset V(G) \times V(G)$
such that if $(v,w) \in E(G)$ then $(w,v) \in E(G)$. Hence our
graphs are undirected and do not have multiple edges, but may have
loops (if $(v,v) \in E(G)$). If $(v,w) \in E(G)$ we will often say
that $v$ and $w$ are \textit{adjacent} and denote this as $v \sim
w$. Given a pair of graphs $G$ and $H$, a \textit{graph
homomorphism} (or \textit{graph map}) $f:G \rightarrow H$ is a map
of the vertex set $f:V(G) \rightarrow V(H)$ that preserves
adjacency: if $v \sim w$ in $G$, then $f(v) \sim f(w)$ in $H$.  With
these as our objects and morphisms we obtain a category of graphs
which we will denote ${\mathcal G}$.

If $v$ and $w$ are vertices of a graph $G$, the \textit{distance}
$d(v,w)$ is the length of the shortest path in $G$ from $v$ to $w$.
The $diameter$ of a finite connected graph $G$, denoted $\diam(G)$
is the maximum distance between two vertices of $G$. The
\textit{neighborhood} of a vertex $v$, denoted $N_G(v)$ (or $N(v)$
if the context is clear), is the set of vertices that are precisely
distance 1 from $v$ (so that $v \in N(v)$ if and only if $v$ has a
loop).  If $v$ and $w$ are vertices of a graph such that $N(v)
\subseteq N(w)$ then we call the map $f:G \rightarrow G \backslash
v$ that sends $v$ to $w$ a \textit{folding} of the vertex $v$; we
will also say that $G$ \textit{folds onto} the graph $G \backslash
v$.

There are several simplicial complexes one can associate with a
given graph $G$.  One such construction is the \textit{clique
complex} $\Delta(G)$, a simplicial complex with vertices given by
all \textit{looped} vertices of $G$, and with faces given by all
cliques (complete subgraphs) on the looped vertices of $G$.

We next recall the definition of our main object of study, the
$\Hom$ complex. (Versions of) this construction were originally used
by Lov{\'a}sz, Babson and Kozlov, and others to provide so-called
\textit{topological} lower bounds to the chromatic numbers of graphs
(see ~\cite{Kchr} for a nice survey). We will use the following
definition.

\begin{defn}
For graphs $G$ and $H$, we define $\Hom(G,H)$ to be the poset whose
elements are given by all functions $\eta: V(G) \rightarrow 2^{V(H)}
\backslash \{\emptyset\}$, such that if $(x,y) \in E(G)$, then for
any $\tilde x \in \eta(x)$ and $\tilde y \in \eta(y)$ we have
$(\tilde x, \tilde y) \in E(H)$.  The partial order is given by
containment, so that $\eta \leq \eta^\prime$ if $\eta(x) \subseteq
\eta^\prime(x)$ for all $x \in V$.
\end{defn}

\begin{center}

\epsfig{file=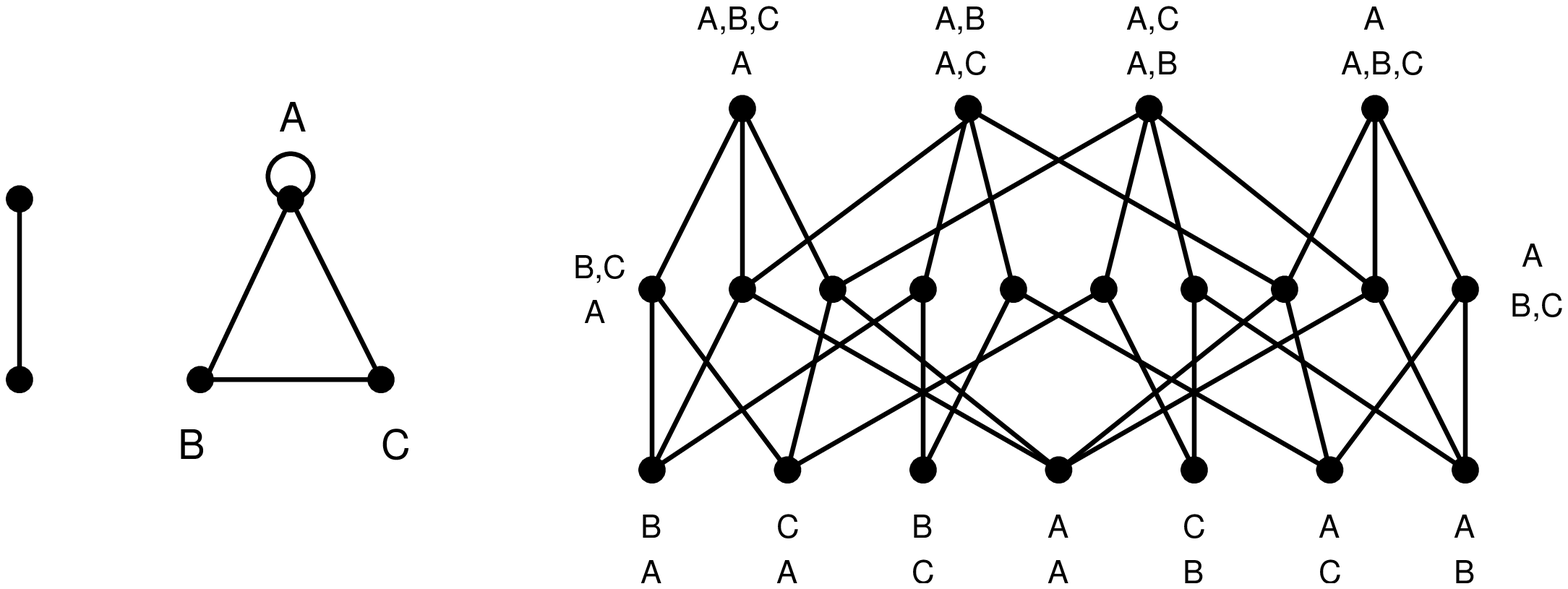, height=1.5 in, width = 5 in}

{The graphs $G$ and $H$, and the poset $\Hom(G,H)$}
\end{center}

One can check that for a fixed graph $G$, $\Hom(G,?)$ (resp.
$\Hom(?,H)$) is a covariant (resp. contravariant) functor from the
category of graphs to the category of posets. We will often speak of
topological properties of the $\Hom$ \textit{complex}.  In this
context we will mean the space obtained as the geometric realization
of the \textit{order complex} of the poset $\Hom(G,H)$.  When the
context is clear, we will refer to this topological space
(realization of a simplicial complex) with the same $\Hom(G,H)$
notation.

\begin{center}

\epsfig{file=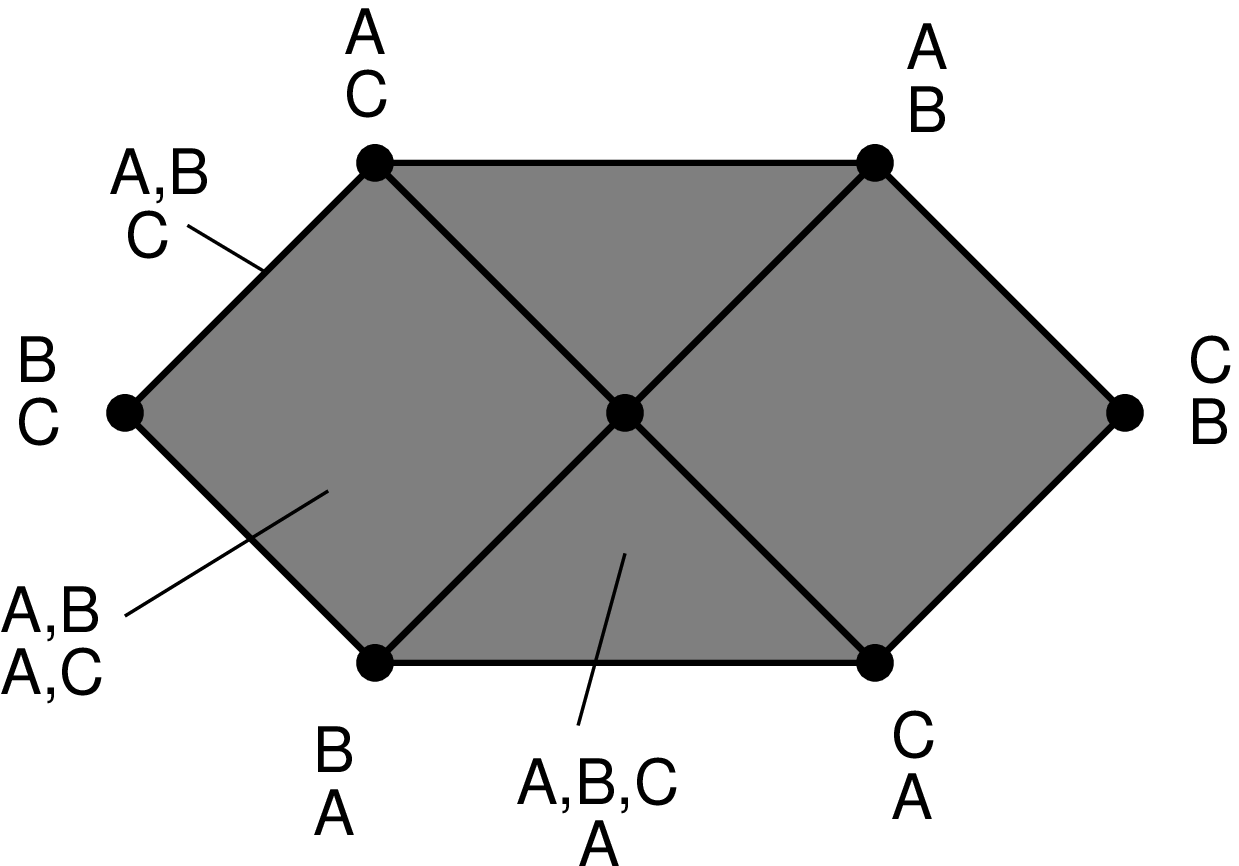, height=1.5 in, width = 2 in}

{The realization of the poset $\Hom(G,H)$ (up to barycentric
subdivision)}
\end{center}

The category ${\mathcal G}$ has a product with a right adjoint given
by the exponential graph construction. We recall these constructions
below.

\begin{defn}
If $G$ and $H$ are graphs, then the \textit{categorical product} $G
\times H$ is a graph with vertex set $V(G) \times V(H)$ and
adjacency given by $(g,h) \sim (g^\prime, h^\prime)$ in $G \times H$
if both $g \sim g^\prime$ in $G$ and $h \sim h^\prime$ in $H$.
\end{defn}

\begin{defn}
For graphs $G$ and $H$, the \textit{categorical exponential graph}
$H^G$ is a graph with vertex set $\{f:V(G) \rightarrow V(H)\}$, the
collection of all vertex set maps, with adjacency given by $f \sim
f^\prime$ if whenever $v \sim v^\prime$ in $G$ we have $f(v) \sim
f^\prime(v^\prime)$ in $H$.
\end{defn}

The exponential graph construction provides a right adjoint to the
categorical product.  This gives the category of graphs the
structure of an \textit{internal hom} associated with the (monoidal)
categorical product (see ~\cite{Doc} for the meaning of these
statements). It turns out that the $\Hom$ complex interacts well
with this adjunction, as described in the following proposition (see
~\cite{Kchr} or ~\cite{Doc} for a proof).

\begin{prop} \label{adjoint}
For $A, B, C$ any graphs, $\Hom(A \times B, C)$ can be included in
$\Hom(A, C^B)$ so that $\Hom(A \times B, C)$ is a strong deformation
retract of $\Hom(A, C^B)$. In particular, we have $\Hom(A \times B,
C) \simeq \Hom(A, C^B)$.
\end{prop}

Note that, as a result of the proposition, we have $\Hom(G,H) \simeq
\Hom({\bf 1}, H^G)$, where ${\bf 1}$ is the graph with a single
looped vertex.  The latter space is homeomorphic to the
(realization) of the clique complex $\Delta(H^G)$.  Hence, up to
homotopy type, the space $\Hom(G,H)$ is just the clique complex on
the looped vertices of the graph $H^G$.  We will use this
identification in the proof of the main theorem.

\section{Proof of the main theorem}

In this section we provide the proof of Theorem ~\ref{maintheorem}.
Note that if $T = {\bf 1}$ is a single looped vertex, we have
$\Hom({\bf 1}, G) \simeq \Delta(G)$, the clique complex on the
looped vertices of $G$.  Hence to obtain the result in this case, we
define $G_{k,X}$ to be the graph obtained by taking the 1-skeleton
of $bd(X) = X^1$, the first barycentric subdivision of the given
complex $X$.  Since the barycentric subdivision of a simplicial
complex is a flag complex, we get that the 1-skeleton provides an
inverse to the $\Delta$ functor in this case, and hence $X \simeq
X^1 = \Delta(G_{k,X})$.

In the general case we will similarly obtain $G_{k,X}$ as the looped
1-skeleton of some iterated subdivision of $X$, but this time we
have to take into account the diameter of the test graph $T$. Recall
the setup: we are given a connected graph $T$ with at least one
edge, and a finite simplicial complex $X$.  If $d = \diam(T)$ is the
diameter of $T$, we fix an integer $k \geq 2$ such that

\begin{center}
$2^{k-1} - 1 \geq d$.
\end{center}

Next, we let $X^k = \bd^k(X)$ denote the $k^{th}$ barycentric
subdivision of the simplicial complex $X$.  We define $G_{k,X}$ to
be the graph given by the 1-skeleton of $X^k$, with loops placed at
every vertex.

\begin{center}

\epsfig{file=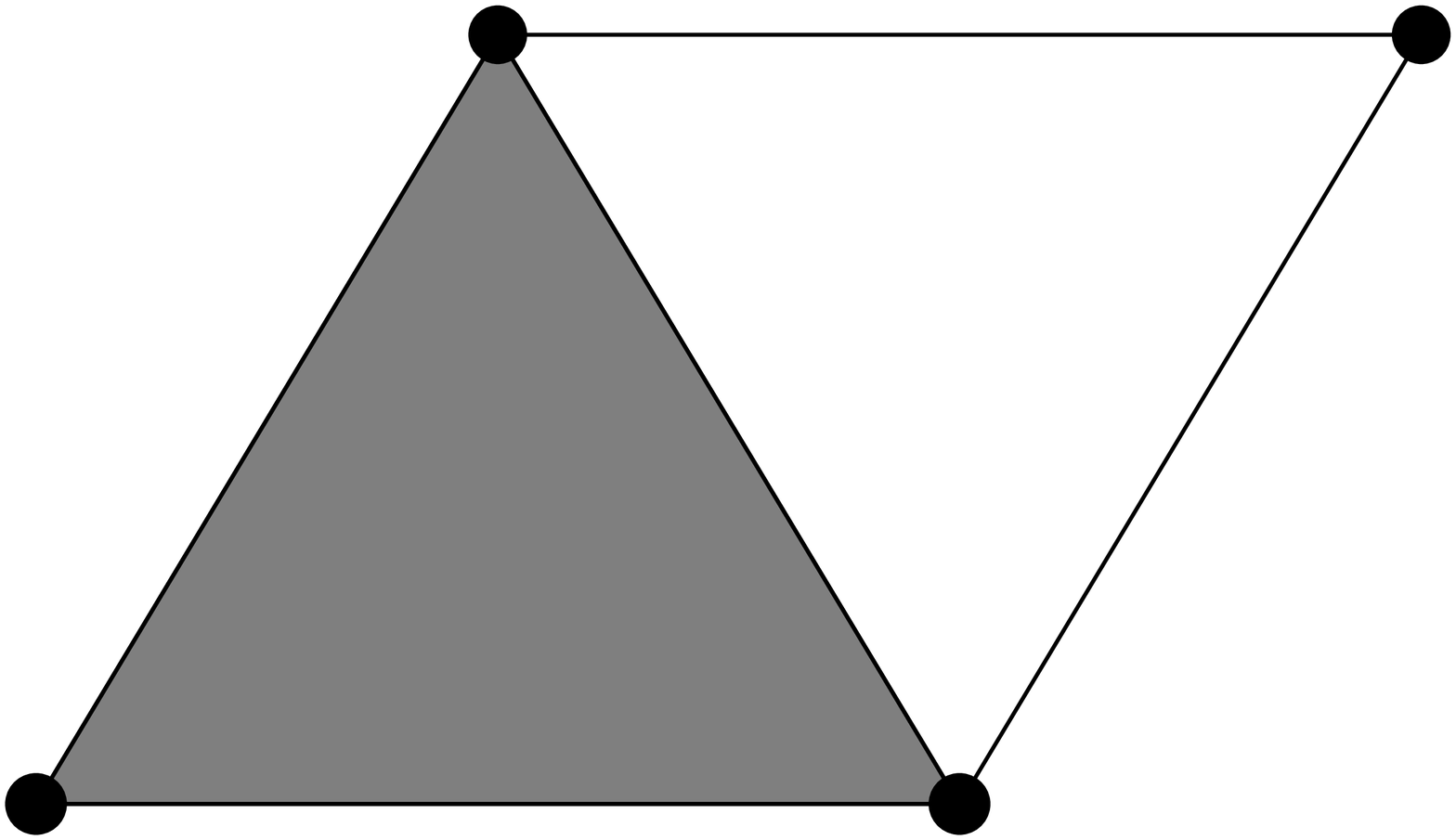, height=1 in, width = 1.75 in}
\epsfig{file=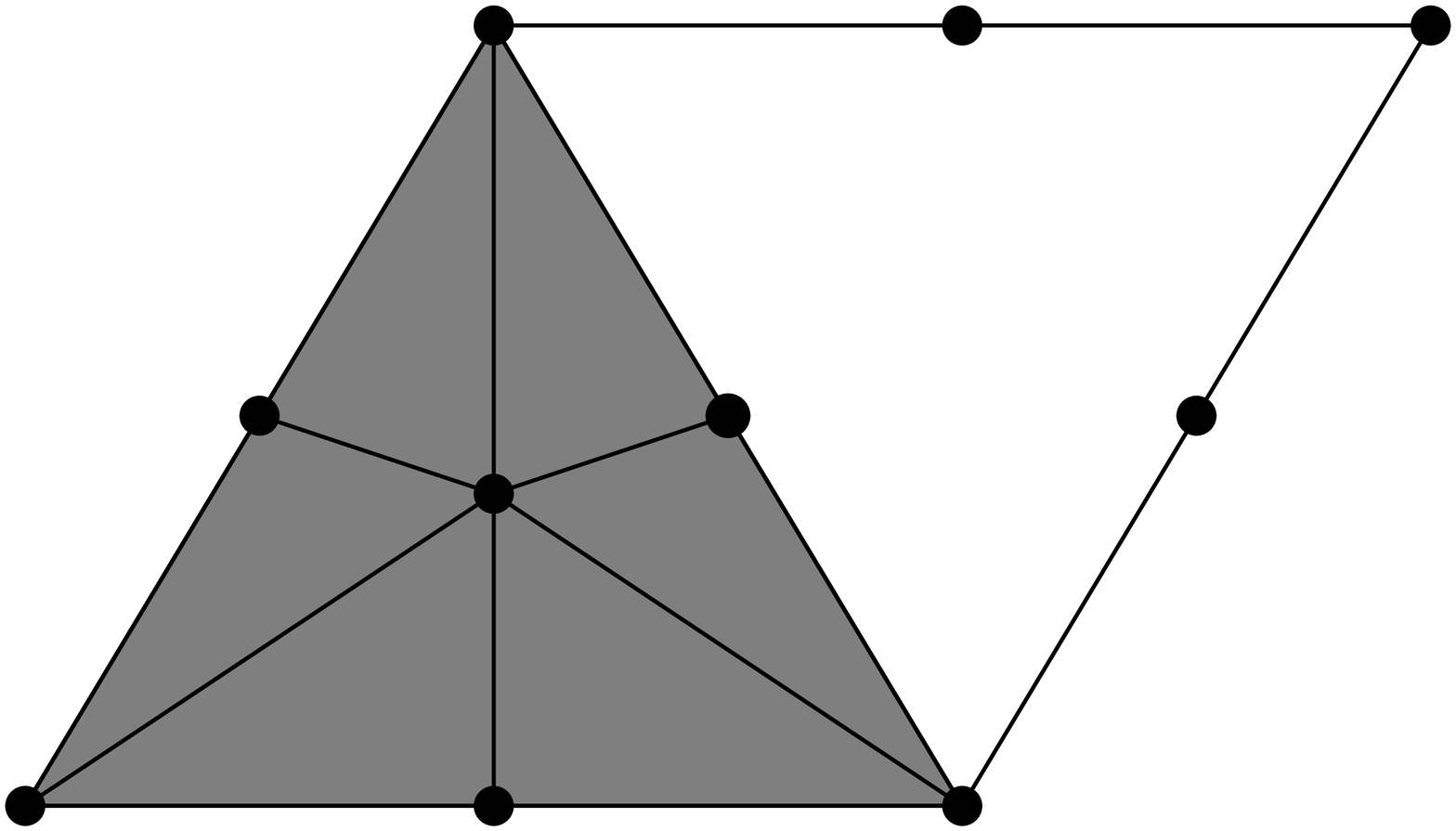, height=1 in, width = 1.75 in}

\end{center}

\vspace{.2 in}

\begin{center}
\epsfig{file=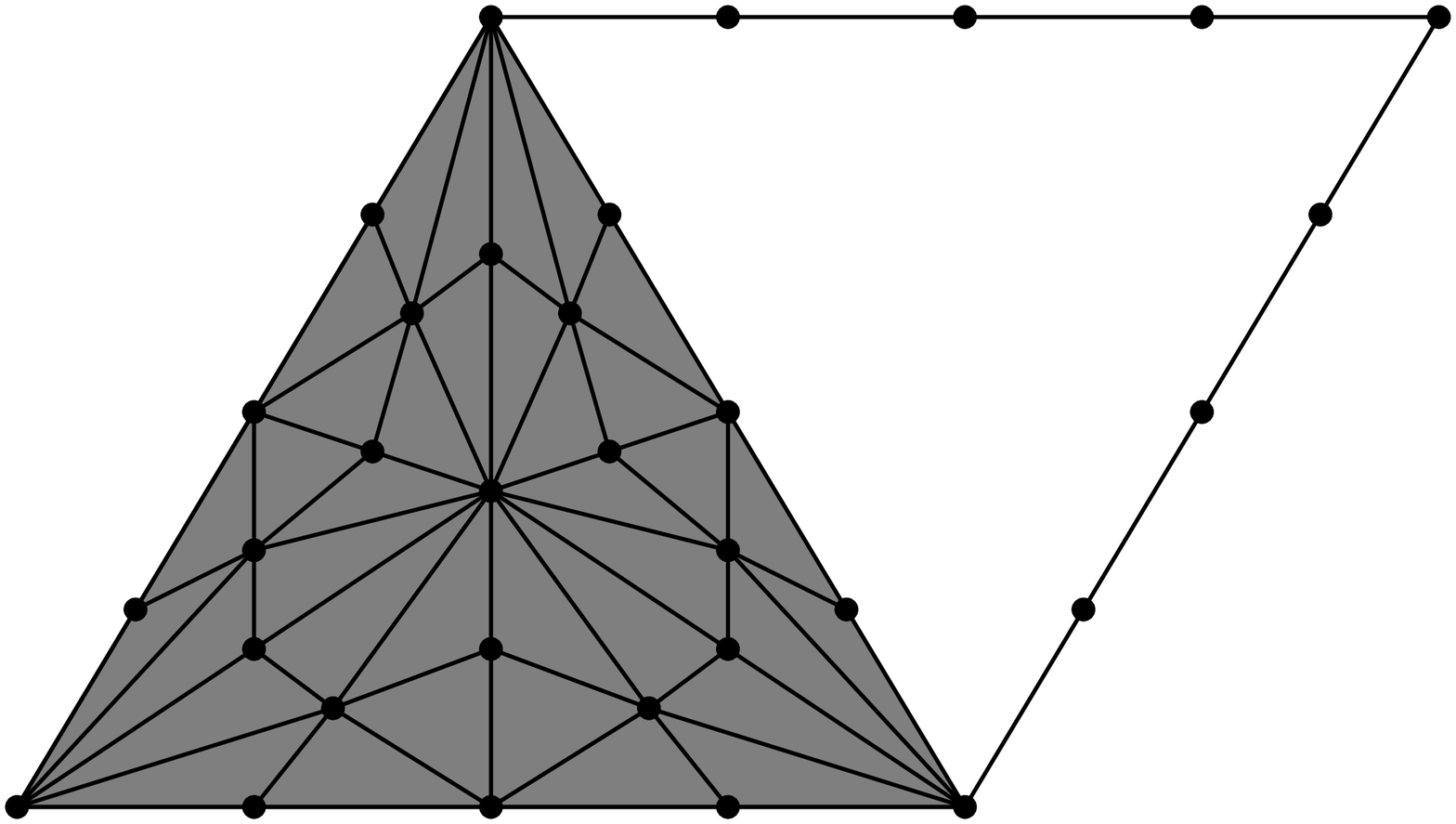, height=1 in, width = 1.75 in}
\epsfig{file=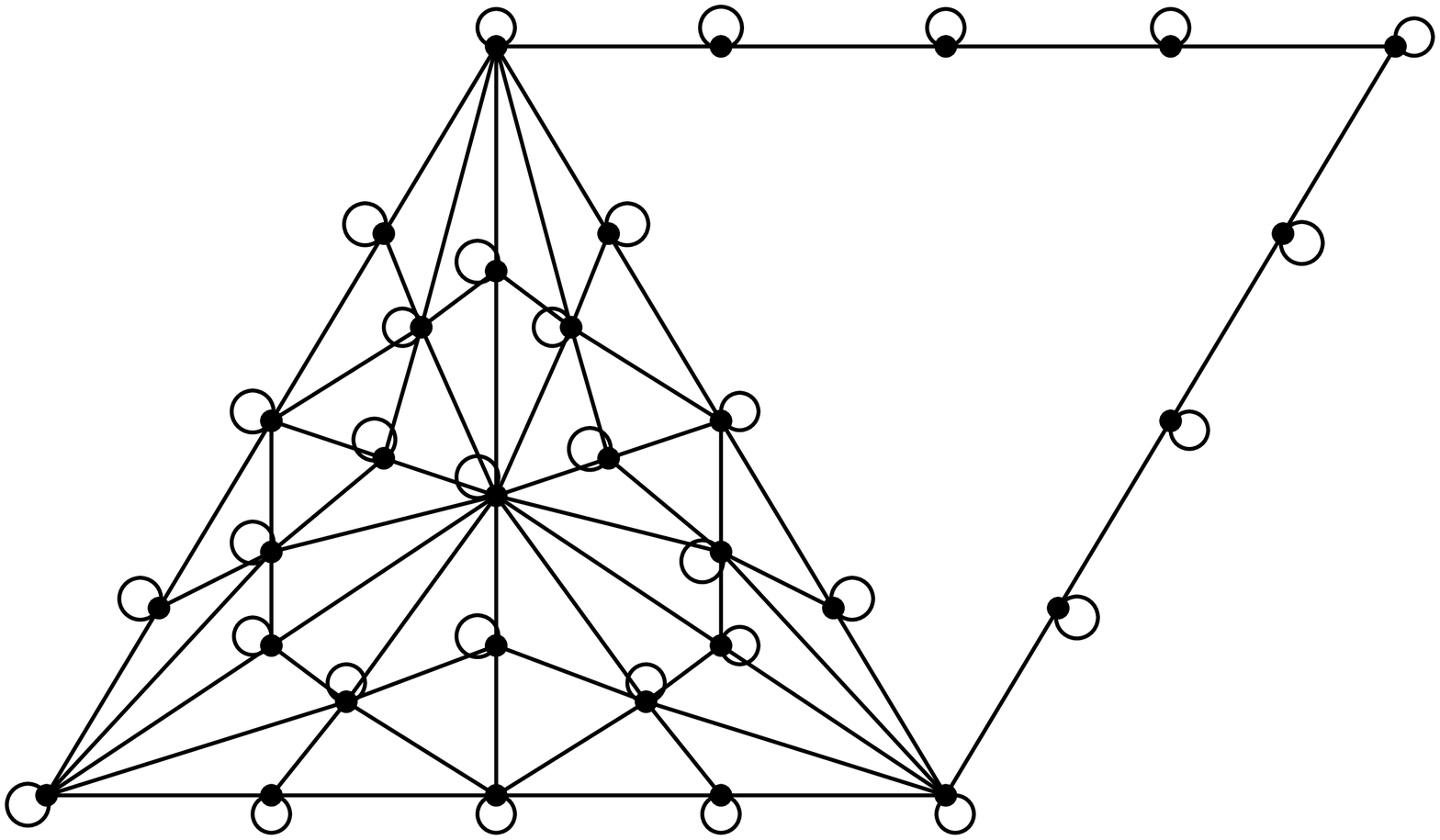, height=1 in, width = 1.75 in}

{The complexes $X$, $X^1$, $X^2$, and the reflexive graph
$G_{2,X}$.}
\end{center}

We claim that $\Hom(T,G_{k,X}) \simeq X$.  From Proposition
~\ref{adjoint}, we have $\Hom(T,G_{k,X}) \simeq \Hom\big({\bf 1},
(G_{k,X})^T\big)$ (where ${\bf 1}$ is the graph with one looped
vertex). The latter space is homeomorphic to
$\Delta\big((G_{k,X})^T\big)$, the clique complex on the (looped
vertices of the) graph $\Delta\big((G_{k,X})^T\big)$. Hence to prove
the main result (Theorem \ref{maintheorem}) it is enough to prove
the following restatement.

\begin{thm}
Let $T$ be an arbitrary connected graph with at least one edge, and
let $X$ be a finite simplicial complex.  Then for $k \geq \Max\{2,
\diam(T)\}$ there is a homotopy equivalence

\begin{center}
$X \simeq \Delta\big((G_{k,X})^T\big)$.
\end{center}

\end{thm}

\begin{proof}
We consider subcomplexes of $\Delta\big((G_{k,X})^T\big)$ of the
form $\Delta\big((G_{k,X}^x)^T\big)$ (see Definition
~\ref{subcomplex} below for the definition of the graph
$G_{k,X}^x$). By Lemma ~\ref{cover} the collection of these
subcomplexes form a cover of $\Delta\big((G_{k,X})^T\big)$, and by
Lemma ~\ref{nerve} the nerve of this cover is isomorphic to the
simplicial complex $X$. By Lemma ~\ref{subcomplex} and Lemma
~\ref{intersection}, these subcomplexes and all nonempty
intersections are contractible.  The result follows from the nerve
lemma of ~\cite{Bjo95}.
\end{proof}

We next turn to the definition of our subcomplexes and the proofs of
the lemmas mentioned above.  Recall that the simplicial complex
$\Delta\big((G_{k,X})^T\big)$ is determined by its 1-skeleton
$(G_{k,X})^T$, whose vertices are given by all graph maps $f:T
\rightarrow G_{k,X}$, and with edges $\{f, f^\prime\}$ whenever
$f(t) \sim f^\prime(t^\prime)$ for all $t \sim t^\prime$ in $T$. We
note that the vertices of the original complex $X$ are naturally
vertices of the graph $G_{k,X}$.  We will work with certain graph
theoretic `open neighborhoods' of these vertices, as described in
the following definition.

\begin{defn} \label{subcomplex}
For a fixed vertex $x$ of the original complex $X$, define
$G^x_{k,X}$ to be the subgraph of $G_{k,X}$ induced by the vertices
$\{w \in G_{k,X}: d(x,w) \leq 2^k-1\}$.

\noindent
Hence the vertices of $G^x_{k,X}$ are the vertices of
$G_{k,X}$ that are distance at most $2^k - 1$ from the vertex $x$.

\begin{center}

\epsfig{file=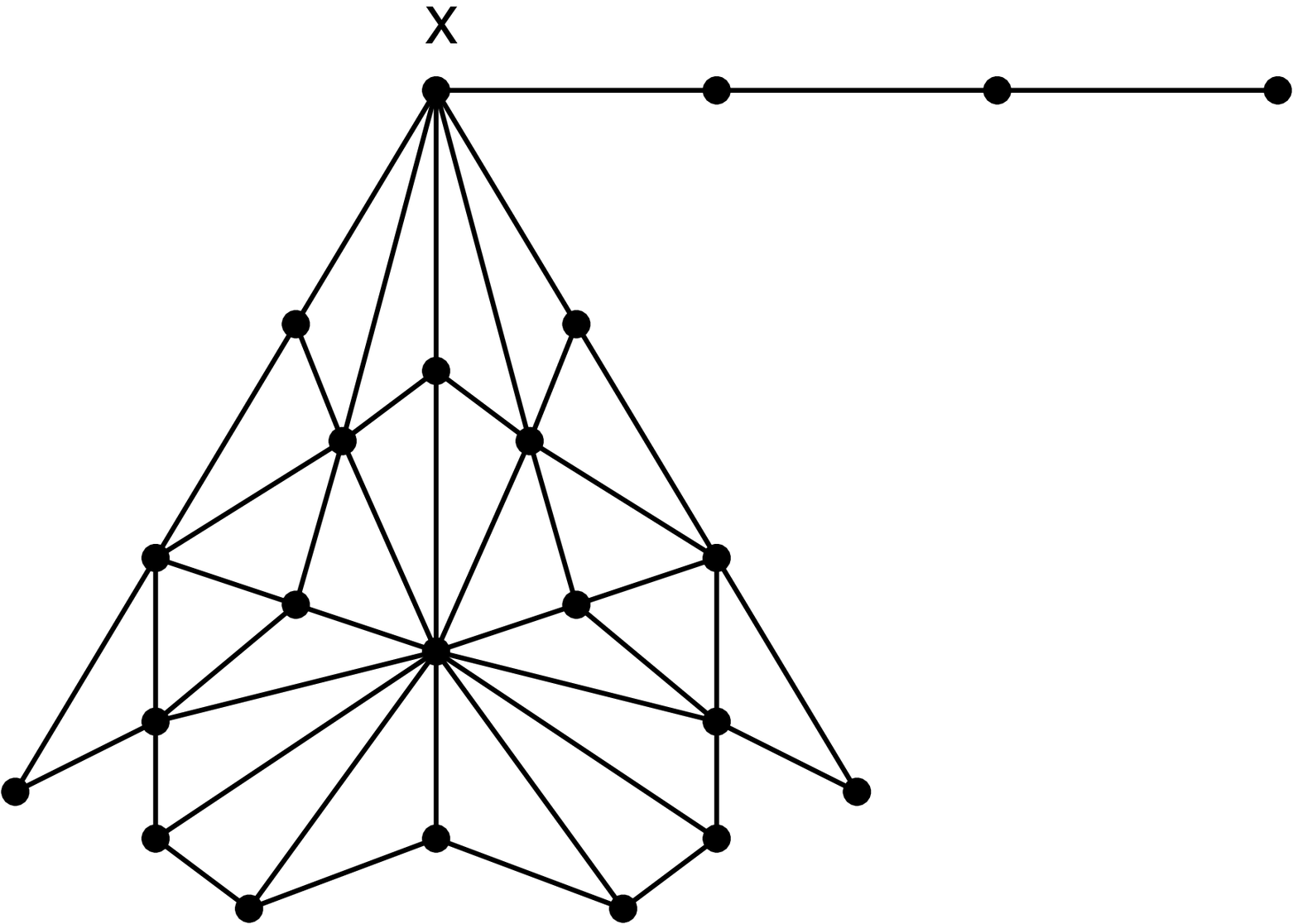, height=1.5 in, width = 2.5 in}

{The graph $G^x_{2,X}$ (without the loops)}

\end{center}

\end{defn}

It is this collection of subcomplexes
$\Big\{\Delta\big((G_{k,X}^x)^T\big)\Big\}_{x \in V(X)}$ that we
wish to show cover the complex $\Delta\big((G_{k,X})^T\big)$.  For
this we will need a general lemma regarding clique complexes of
exponential graphs. For graphs $T$ and $G$, and a simplex $\alpha =
\{f_1, \dots, f_a \} \in \Delta(G^T)$, define $G_\alpha$ to be the
subgraph of $G$ induced by the vertices $\{f_i(t):1 \leq i \leq a, t
\in V(T)\}$. We then make the following observation.

\begin{lemma} \label{diameter}
Let $T$ be a finite connected graph with diameter $d = \diam(T)$,
and suppose $G$ is any graph.  Then $\diam(G_\alpha) \leq \Max \{2,
d\}$ for all $\alpha \in \Delta(G^T)$.
\end{lemma}

\begin{proof}
Suppose $T$ and $G$ are as above, and suppose $\alpha= \{f_1, \dots
f_a\}$ is a face of $\Delta(G^T)$.  Let $v = f_i(t)$ and $v^\prime =
f_{i^\prime}(t^\prime)$ be any two elements of $G_\alpha$. We will
find a path in $G_\alpha$ from $v$ to $v^\prime$ of length $\leq d$.
If $t \neq t^\prime$, then by assumption we have a path in $T$ from
$t$ to $t^\prime$ given by $(t = t_0, t_1, \dots, t_j = t^\prime)$,
with $j \leq d$.  If $t = t^\prime$, we take our path to be $(t,
t_1, t_2 = t)$, where $t_1$ is any neighbor of $t$.  So we have $j
\leq \Max \{2, d\}$

Now, since $\alpha$ is a clique in the graph $G^T$, we have that
$f_i \sim f_j$ for all $1 \leq i,j \leq k$, and hence $f_i(t) \sim
f_j(t^\prime)$ for all adjacent $t \sim t^\prime$.  Hence we can
take our desired path to be $f_i(t) = f_i(t_0), f_i(t_1), \dots,
f_i(t_{j-1}), f_{i^\prime}(t_j) = f_{i^\prime}(t^\prime)$.
\end{proof}

We can now show that our subcomplexes indeed form a cover.

\begin{lemma} \label{cover}
The collection of complexes
$\Big\{\Delta\big((G^x_{k,X})^T\big)\Big\}_{x \in V(X)}$ covers the
complex $\Delta\big((G_{k,X})^T\big)$.
\end{lemma}

\begin{proof}
To simplify indices, in our notation for graphs we will suppress
reference to the integer $k$ and the simplicial complex $X$, so that
for this proof $G = G_{k,X}$ and $G^x = G_{k,X}^x$. If $\alpha$ is a
face of $\Delta(G^T)$ then by Lemma ~\ref{diameter} we have either
$k = 2$ and $\diam(G_\alpha) = 2$, or else $diam(G_\alpha) \leq d
\leq 2^{k-1}-1$. We claim that $G_\alpha \subset G^x$ for some $x
\in V(X)$, which would prove our claim.

Let $m = \Min\{d(w,x): w \in G_\alpha, x \in V(X)\}$. Note that $m
\leq 2^{k-1}$ since every vertex of $X^k$ is within distance
$2^{k-1}$ of some vertex of the original complex $X$.

If $m = 0$ then we have $y \in G_\alpha$ for some vertex $y \in
V(X)$. Hence $G_\alpha \subset G^y$ since $G^y$ contains all
vertices distance at most $2^k - 1$ from $y$ (this number is at
least 2 since $k \geq 2$).

If $m > 0$ let $w$ be a vertex of $G_\alpha$ such that $d(w,x) = m$
for some vertex $x \in X$, and choose $w$ such that it is contained
in the interior of a face of $X$ of minimum dimension.  We need to
show that $G_\alpha \subset G^x$.  To see this, first consider the
case that $k > 2$.  By Lemma ~\ref{diameter}, all vertices
$w^\prime$ in $G_\alpha$ are distance at most $d \leq 2^{k-1} - 1$
from $w$. So all vertices of $G_\alpha$ are at most $m + d \leq
2^{k-1} + 2^{k-1} - 1 = 2^k - 1$ away from $x$, which implies
$G_\alpha \subset G^x$.

If $k = 2$, then we have $m=1$ or $m=2$.  If $m = 1$ then all
vertices of $G_\alpha$ are distance at most $1 + 2 = 3 = 2^k - 1$
away from $x$, as desired.  If $m = 2$, then all vertices of
$G_\alpha$ are distance at least 2 from \textit{every} vertex of
$X$.  Now, $w$ is contained in the interior of some face $F_w = \{x,
x_1, \dots, x_j\}$ of the original complex $X$.   If $w^\prime$ is
any other vertex of $G_\alpha$, then $w^\prime$ cannot be contained
in any proper face of $F_w$ since otherwise we would have taken $w =
w^\prime$.  Hence $w^\prime$ is contained in the interior of $F_w$,
so that $d(w^\prime, x) \leq 2^k - 1$, as desired.  This shows that
$G_\alpha \subset G^x$.

\end{proof}

We next turn to the combinatorics of this cover.  Recall that the
\textit{nerve} of a covering by subcomplexes is the simplicial
complex with vertices given by the subcomplexes and with faces
corresponding to all non-empty intersections.  We then have the
following observation.

\begin{lemma} \label{nerve}
The nerve of the covering of $\Delta\big((G_{k,X})^T\big)$ given by
the subcomplexes $\Delta\big((G_{k,X}^x)^T\big)$ is isomorphic to
the simplicial complex $X$.
\end{lemma}

\begin{proof}
By construction, the vertices of the nerve determined by the
$\Delta\big((G_{k,X}^x)^T\big)$ are indexed by $V(X)$, the vertices
of the simplicial complex $X$.  A collection $I \subseteq V(X)$ of
such subcomplexes has nonempty intersection if and only if there
exists a vertex $x$ within distance $2^k - 1$ from each $v \in I$ in
$X^k$, the $k^{th}$ barycentric subdivision of $X$.  But this occurs
if and only if the collection $I$ of vertices form a face of $X$.
\end{proof}

Next we wish to show that each subcomplex
$\Delta\big((G_{k,X}^x)^T\big)$ is contractible. To do this we will
show that each graph $G_{k,X}^x$ is in fact \textit{dismantlable}.
 Recall that a finite graph $G$ is called \textit{dismantlable} if it
can be folded down to the looped vertex ${\bf 1}$ (see ~\cite{BW04}
and ~\cite{Doc} for other characterizations). It follows from the
results of ~\cite{Kcol} that if $G$ is dismantlable, then
$\Hom(S,G)$ is contractible for \textit{any} graph $S$. Hence to
show that the subcomplexes $\Hom(T,G_{k,X}^x) \simeq
\Delta\big((G_{k,X}^x)^T\big)$ are each contractible, it suffices to
show that each graph $G_{k,X}^x$ is dismantlable.

For this we will describe a recursive folding procedure for the
graph $G_{k,X}^x$.  In our induction we will need the fact that
barycentric subdivision preserves dismantlability, as described by
the following lemma.

\begin{lemma} \label{subdivision}
If $G$ is a dismantlable graph and $\Delta(G)$ is its clique complex
(on its looped vertices), then the one-skeleton of $\bd(\Delta(G))$
is again dismantlable.
\end{lemma}

\begin{proof}
Suppose $G$ is a dismantlable graph, and let $G^\prime$ denote the
graph obtained by taking the looped one-skeleton of
$\bd(\Delta(G))$.
 We can think of $G^\prime$ as the graph whose vertices are the elements of the poset
$\Hom({\bf 1}, G)$, with adjacency given by $x \sim y$ if $x$ and
$y$ are comparable.

To show that $G^\prime$ is dismantlable, we proceed by induction on
$n$, the number of looped vertices of $G$.  If $n = 1$ we have that
$G = G^\prime$ is a single looped vertex, and hence dismantlable.

Next suppose $n > 1$, and let $v$ and $w$ be distinct looped
vertices of $G$ such that $N_G(v) \subseteq N_G(w)$.  For future
reference, we let $N_G(v) = \{v, w, v_1, \dots, v_m \}$ denote the
neighboring vertices of $v$ in the graph $G$.  We will use the
following running example, in which the loops (present on all
vertices) will be omitted for the sake of space.

\begin{center}

\epsfig{file=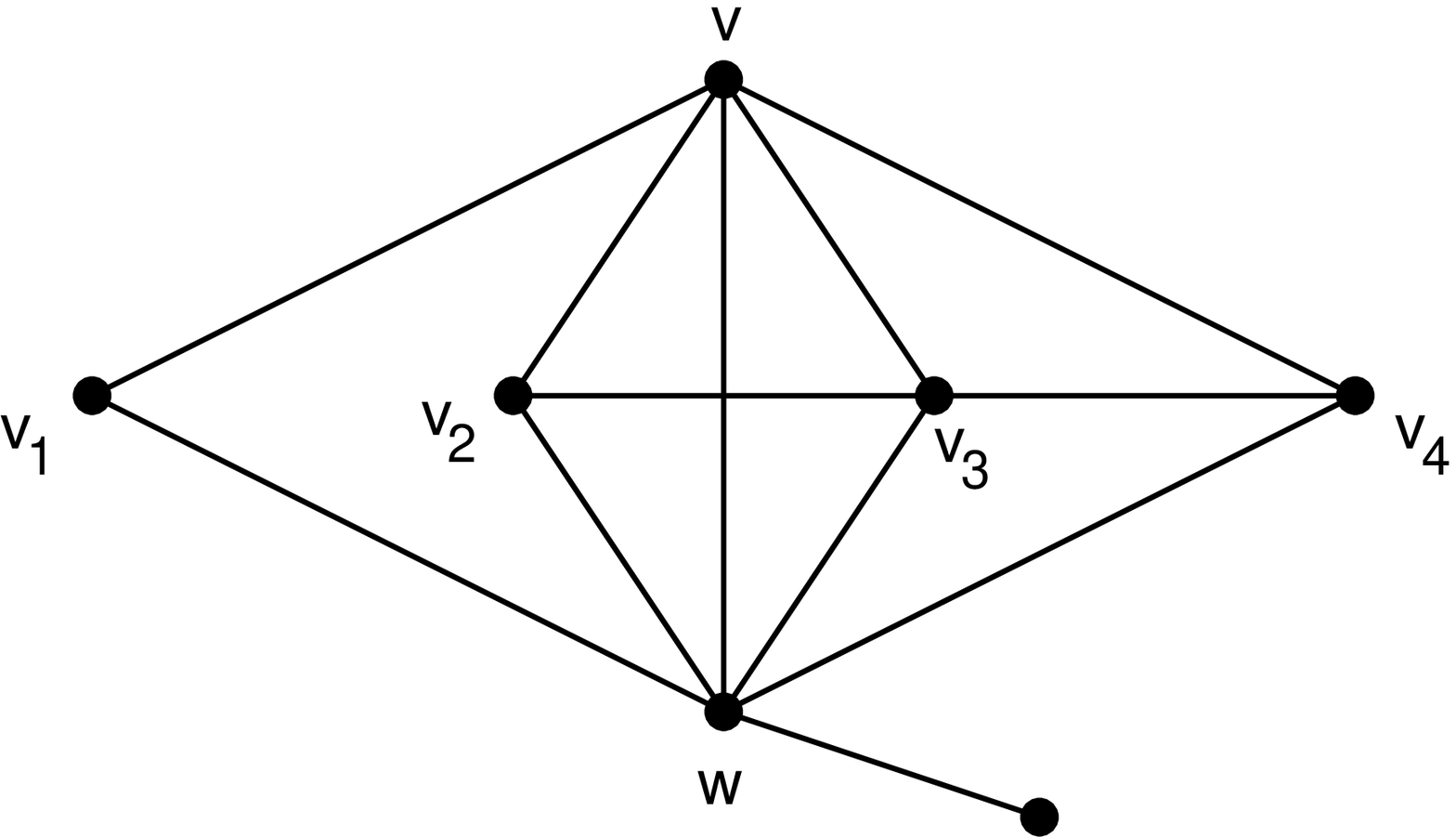, height=1.25 in, width = 2 in}

{The containment $N_G(v) \subseteq N_G(w)$.}
\end{center}

For the inductive step, we need to fold away all vertices in $G^\prime$ that are
barycenters of simplices that have $v$ as a vertex (including the
vertex $v$ itself).  But this is precisely $N_{G^\prime}(v)$, the
collection of neighboring vertices of $v$ in the graph $G^\prime$.

We will first fold away the vertices in $N_{G^\prime}(v)$ that are
furthest from $w$.  We let $S$ denote the collection of vertices in
$N_{G^\prime}(v)$ that are barycenters of simplices that do
\textit{not} contain $w$.  So $S$ is the collection of vertices in
$\bd(\Delta(G))$ that are barycenters of simplices with vertices
among the set $\{v, v_1, \dots, v_m\}$.

Each vertex $s \in S$ is the barycenter of a face of a certain
dimension, and we will fold away the elements of $S$ in descending
order according to this dimension.  If $s$ is the barycenter of a
face $\{v, v_{i_1}, \dots, v_{i_r}\}$ of \textit{maximal} dimension
then we have $N_{G^\prime}(s) \subseteq N_{G^\prime}(y)$, where $y
\in G^\prime$ is the barycenter of the face $\{v, v_{i_1}, \dots
v_{i_r}\,w \}$; this collection forms a face of $\Delta(G)$ since
$N_G(v) \subseteq N_G(w)$.  Hence we can fold away $s$ in this case.

In general, $s$ is the barycenter of a face $F_s = \{v, v_{j_1},
\dots, v_{j_\ell}\}$ and, as we have folded away the vertices of
greater dimension in $S$ (barycenters of faces that contain $F_s$),
we have $N(s) \subseteq N(y)$ in the resulting graph, where again
$y$ is the barycenter of the face $\{v, v_{j_1}, \dots, v_{j_\ell},
w\}$.

\begin{center}

\epsfig{file=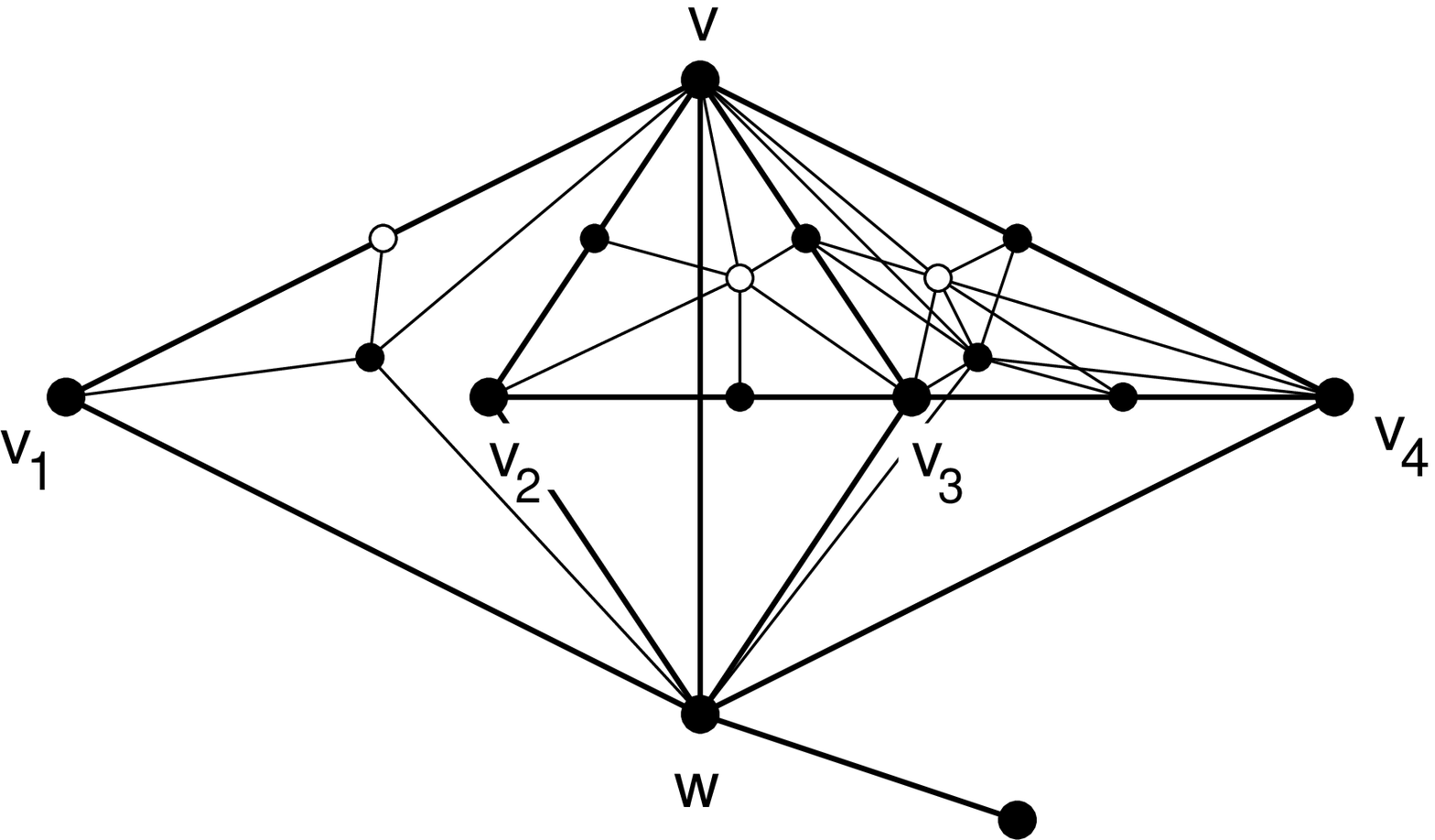, height=1.9 in, width = 2.6 in} \hspace{.1
in} \epsfig{file=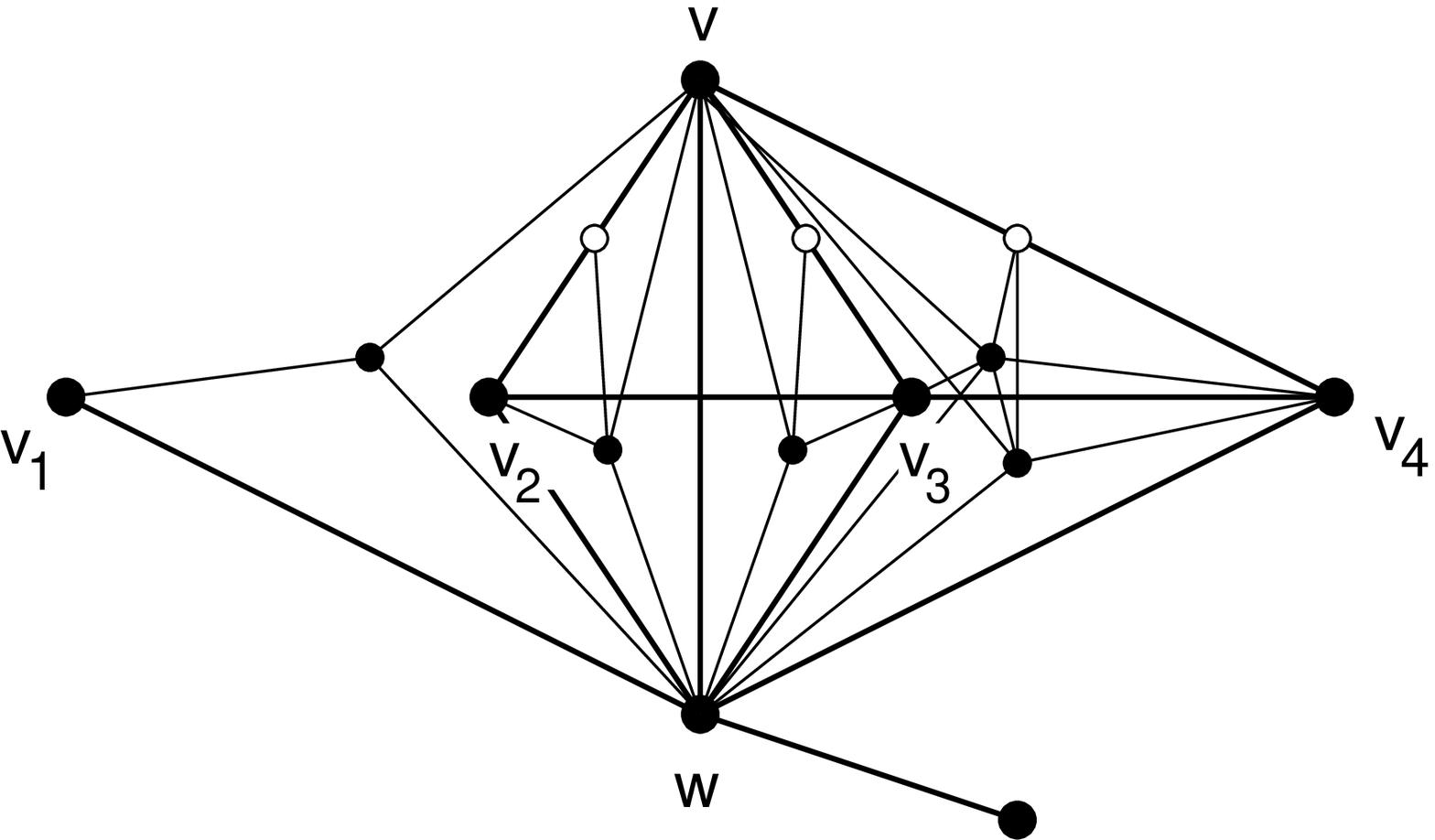, height=1.9 in, width = 2.6 in}

\end{center}

\begin{center}

Folding away the vertices of $S$

\end{center}

In the diagram above, the first step is to fold away the barycenters
of $\{v, v_1\}$, $\{v,v_2,v_3\}$, and $\{v,v_3,v_4\}$ (the vertices
in white). In the second step we fold away the barycenters of
$\{v,v_2\}$, $\{v,v_3\}$, and $\{v,v_4\}$.

Next we fold away the vertex $v$.  If $u \in N(v)$ is a neighbor of
$v$ in the graph at this stage of the folding, then $u$ is the
barycenter of some face that contains both $v$ and $w$, and hence we
have $N(v) \subseteq N(z)$, where $z$ is the barycenter of $\{v,w\}$.  We fold away $v$ and now have that all
neighbors of $z$ are barycenters of faces that contain the vertex
$w$.  Hence we now have $N(z) \subset N(w)$, and we proceed to fold away the
vertex $z$.

\begin{center}

\epsfig{file=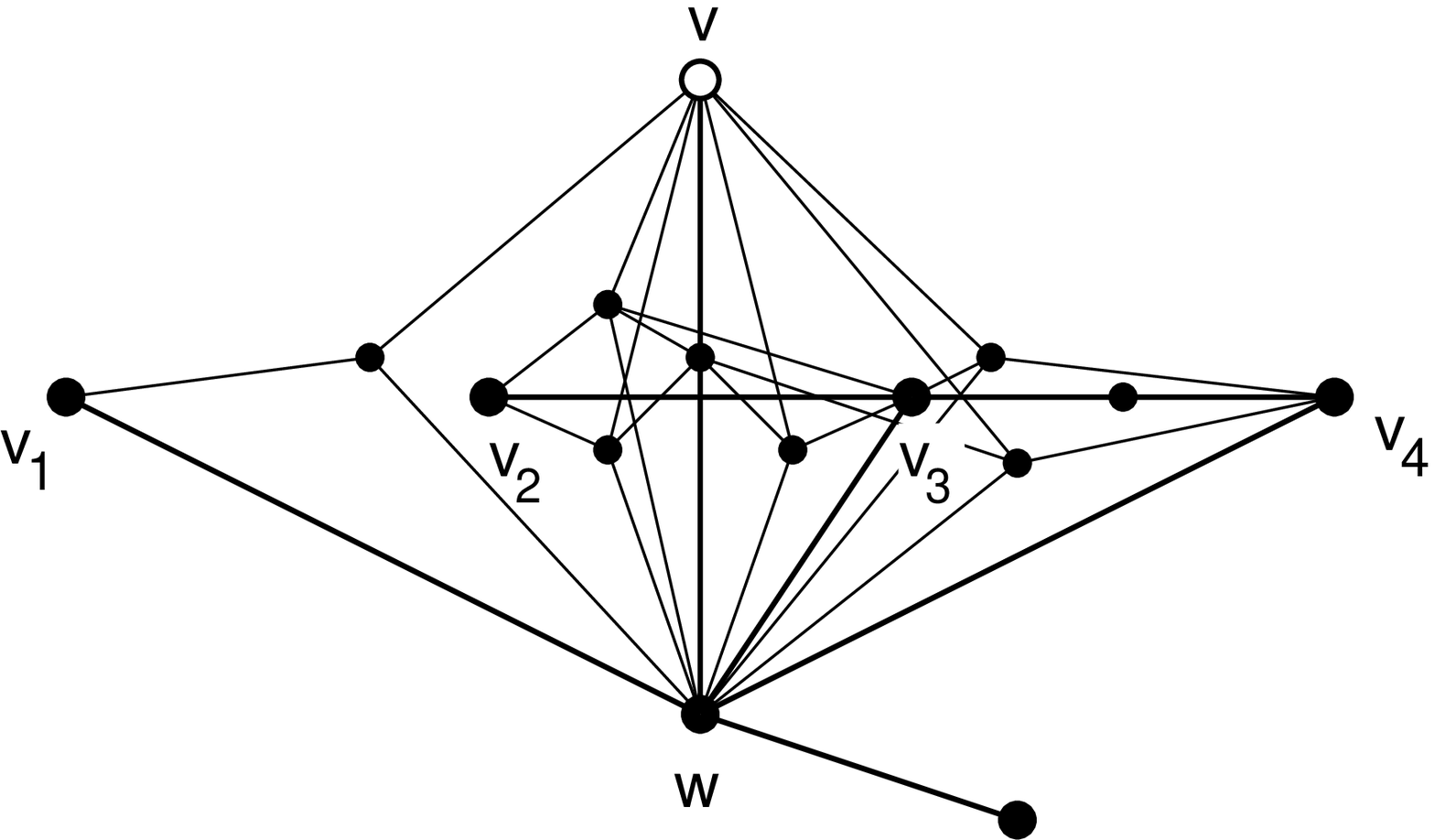, height=1.9 in, width = 2.6 in} \hspace{.1
in} \epsfig{file=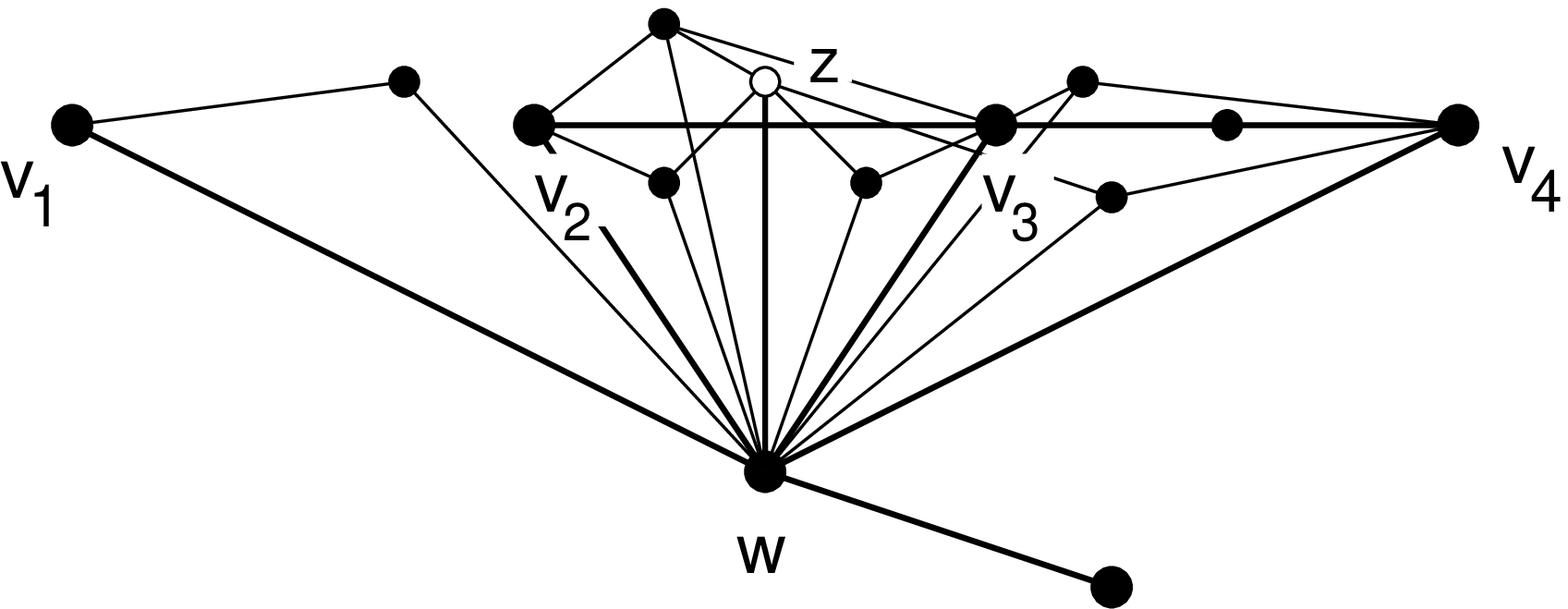, height=1.9 in, width = 2.6 in}

{Folding away $v$ and $z$}
\end{center}

At this point, we are left with a subset $Y \subseteq
N_{G^\prime}(v)$ that consists of vertices that are barycenters of
faces that contain $v$, $w$, and at least one vertex from $\{v_1,
\dots, v_m\}$.   We fold away these vertices in \textit{ascending}
order according to their dimension.  If $y \in Y$ is the barycenter
of a face $\{v, w, v_{i}\}$ of \textit{minimal} dimension, then
$N(y) \subset N(z)$, where $z$ is the barycenter of the face
consisting of $\{w, v_{i}\}$ (since vertices that are barycenters of
faces including $v$ have been folded away). In the general case, $y$
is the barycenter of a face $\{v, w, v_{i_1}, \dots, v_{i_j} \}$
and, as we have folded away the vertices of smaller dimension in
$Y$, we now have $N(y) \subset N(z)$, where again $z$ is the
barycenter of the face $\{w, v_{i_1}, \dots, v_{i_j}\}$.

\begin{center}

\epsfig{file=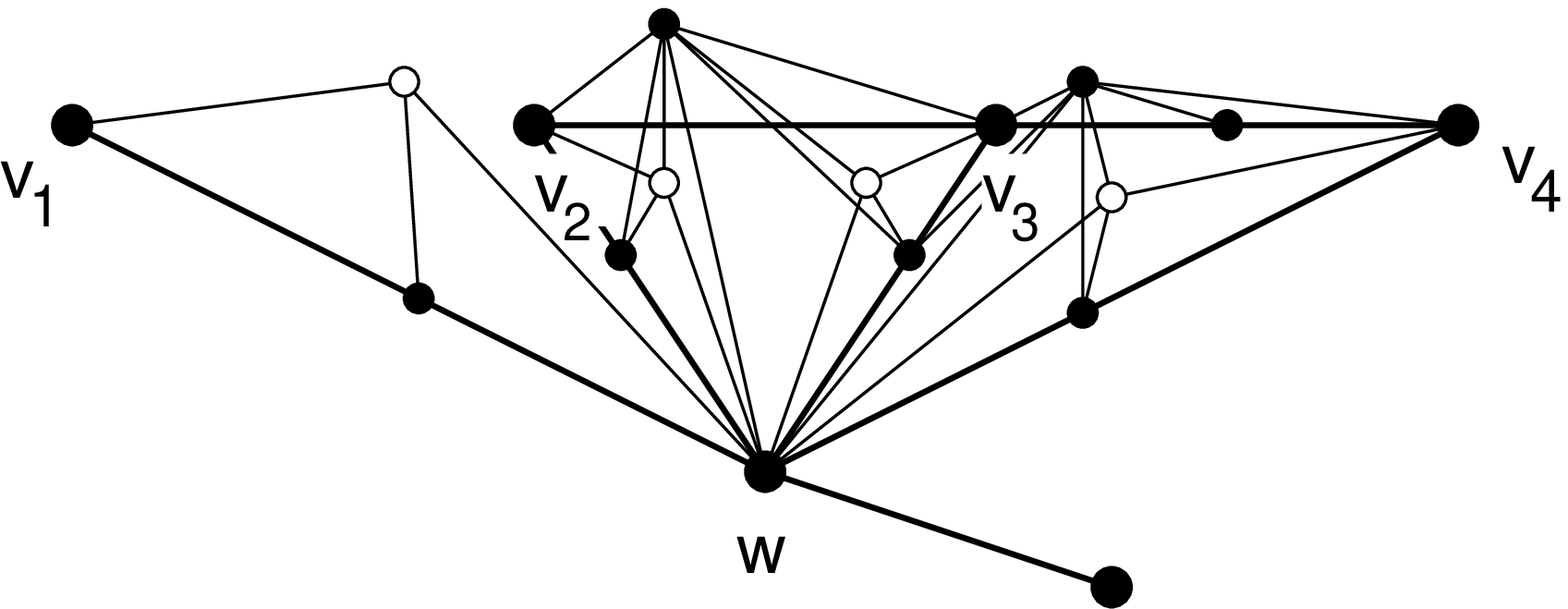, height=1.3 in, width = 2.6 in} \hspace{.1
in} \epsfig{file=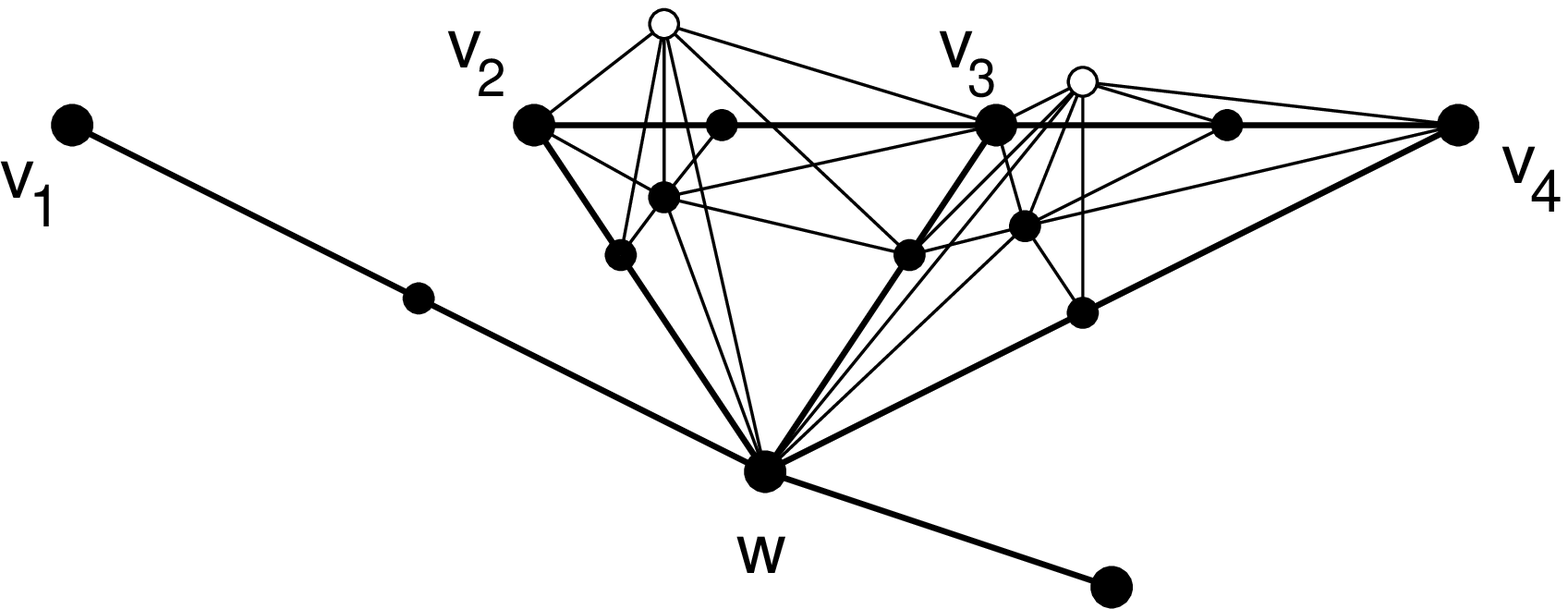, height=1.3 in, width = 2.6 in}

{Folding away the remaining vertices of $N_{G^\prime}(v)$.}
\end{center}

\end{proof}

We can now use this to prove the following result concerning the
$G_{k,X}^x$ graphs.

\begin{lemma} ~\label{subcomplex}
For any vertex $x \in X$, the graph $G_{k,X}^x$ is dismantlable.
\end{lemma}

\begin{proof}
Recall that $G_{k,X}^x$ is the subgraph of $G_{k,X}$ induced by the
vertices that are distance at most $2^k - 1$ from $x$.  We will
prove the claim by induction on $k$.  For $k = 1$ the graph
$G_{k,X}^x$ consists of $N_{G_{k,X}}(x)$, the neighbors of the
vertex $x$ in $G_{k,X}$ (including $x$ itself). Hence $G_{k,X}^x$
folds down to the single looped vertex $x$, as desired.

Next suppose $k > 1$.  Our plan is to first fold away the vertices
in $G_{k,X}^x$ that are distance \textit{exactly} $2^k -1$ from $x$.
The resulting subgraph one obtains is the looped 1-skeleton of the
barycentric subdivision of the clique complex $\Delta(G_{k-1, X}^x)$
(this graph is called $\big(G_{k-1,X}^{x}\big)^\prime$ in the
notation of the proof of Lemma ~\ref{subdivision}). By induction,
together with Lemma ~\ref{subdivision}, this graph is dismantlable
and hence our claim will be proved.

Let $V_x$ denote the collection of vertices in $G_{k,X}^x$ that are
distance exactly $2^k - 1$ from $x$; it is this collection of
vertices that we wish to fold away. First we set up some notation.
Note that every vertex $v$ in the graph $G_{k,X}$ has a pair of
parameters $\alpha(v) = (i,j)$ associated with it, where $i$ is the
dimension of the face in $X$ that $v$ lies in, and where $j$ is the
dimension of the face of $X^{k-1}$ that $v$ is the barycenter of
(note that $j \leq i$).  We will say that $v$ is \textit{of type}
$(i,j)$ if $\alpha(v) = (i,j)$.

\begin{center}

\epsfig{file=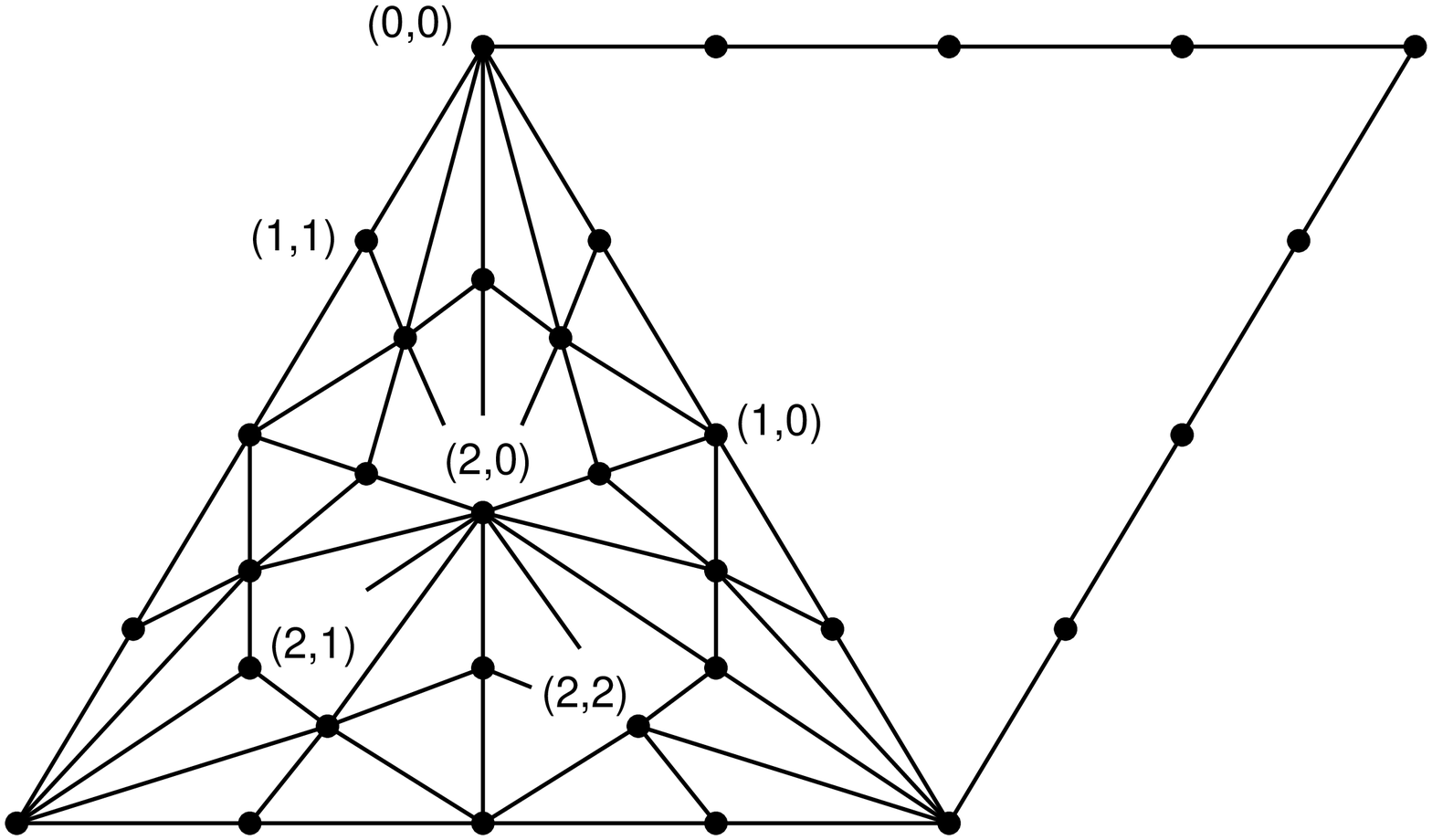, height=1.75 in, width = 3 in}

{The types $(i,j)$ of various vertices in the graph $G_{k,X}$}
\end{center}

We will fold away the vertices of $V_x \subset G_{k,X}^x$ in
lexicographic order according to their type $(i,j)$. First note that
if $v \in V_x$ is of type $(i,j)$ then $j \geq 1$, and hence our
base case to consider is when $v \in V_x$ is of type $(1,1)$.  In
this case $v$ is the barycenter of an edge $\{a,b\}$ in $X^{k-1}$,
where $b$ is a vertex of $X$, and $a$ is distance $2^{k-1}-2$ from
$x$. Any neighbor $w \in N_{G_{k,X}^x}(v)$ of $v$ is a barycenter of
a simplex that has $a$ as a vertex; hence we have $w \sim a$.  We
conclude that $v$ can be folded onto the neighboring vertex $a$.

Next we consider the case $v$ is of type $(i,j)$, where $i > 1$ is
fixed.  We proceed by induction on $j$.  If $j=1$ then $v$ is the
barycenter of an edge $\{c,d\}$, where $c \notin V_x$ and $d$ is of
type $(i,0)$ and is distance $2^k - 2$ from $x$.  Any other neighbor
$w \in N_{G_{k,X}^x}(v)$ is the barycenter of a simplex that has $d$
as a vertex; we conclude that $w \sim d$. Hence in this case $v$ can
be folded onto $d$.

\begin{center}

\epsfig{file=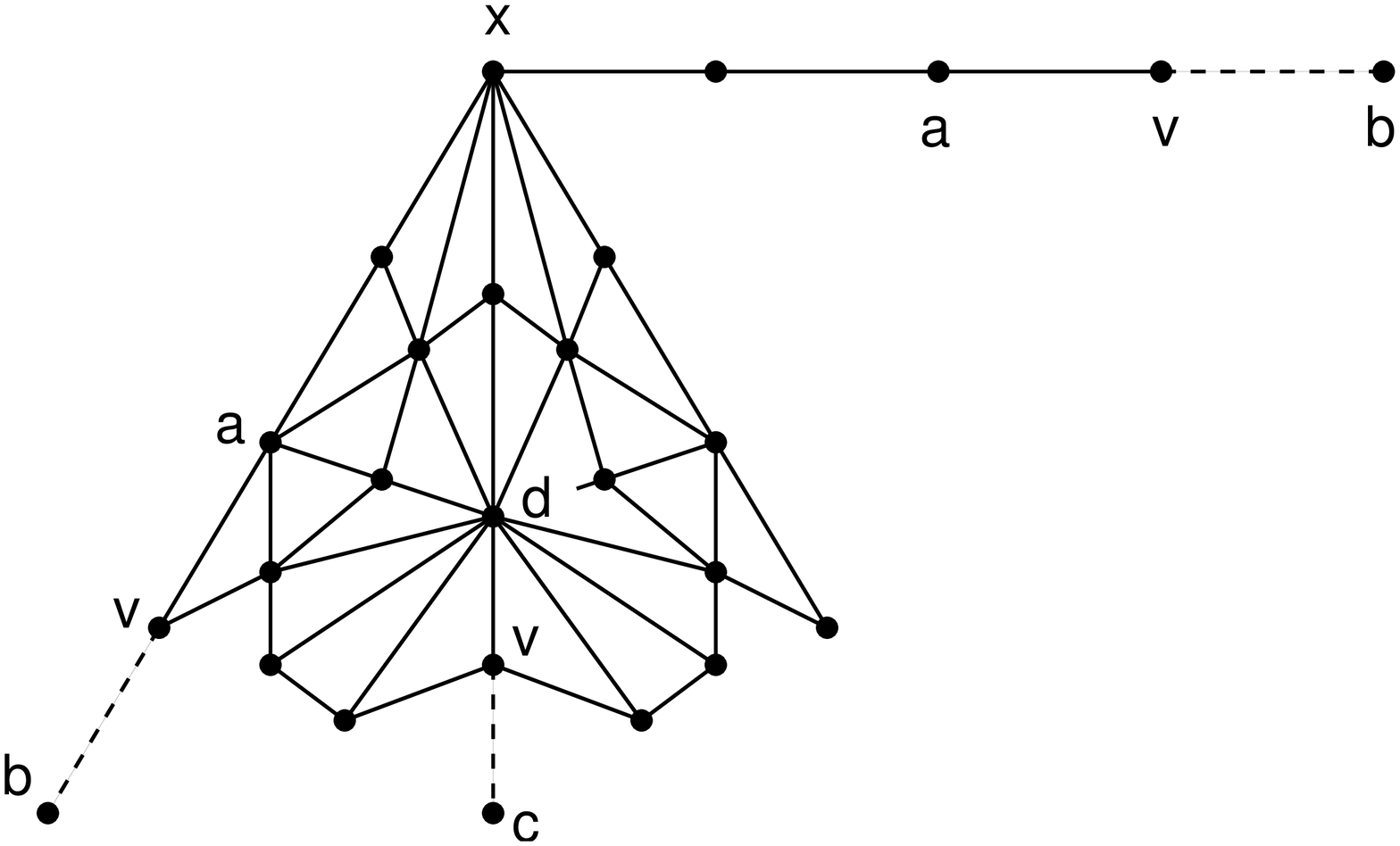, height=1.75 in, width = 3 in}

{Folding away the vertex $v$ when $v$ is of type $(i,1)$}
\end{center}

For the same fixed $i > 1$, we next consider the case that $v$ is of
type $(i,j)$, where $j > 1$.  By induction, we have that all
vertices in $V_x$ of type $(i^\prime, k)$ and of type $(i,j^\prime)$
have been folded away, where $i^\prime < i$ and $j^\prime < j$. Pick
a vertex $w \in N(v)$ in the neighborhood of $v$ such that $w \in
G_{k-1,X}^x$ and such that the type of $w$ is largest in the
lexicographic order - that is, of type $(i,j)$ where $j$ is maximum
among maximum $i$.

We claim that $N_{G_{k,X}^x}(v) \subseteq N_{G_{k,X}^x}(w)$, so that
the vertex $v$ can be folded onto $w$.  To see this, suppose $u \in
N_{G_{k,X}^x}(v)$.  If $u \in V_x$ (so that $d(u,x) = 2^k - 1$),
then by induction we know that $u$ is of type $(i^\prime,
j^\prime)$, where either $i^\prime > i$ or else $i^\prime = i$ and
$j^\prime > j$.  In either case we see that $u$ is the barycenter of
a simplex $U$ that contains the vertex $w$, and hence $u \sim w$ as
claimed. If $u \notin V_x$, so that $d(u,x) = 2^k - 2$, then either
$u = w$ or else the type of $u$ is lexicographically smaller than
the type of $w$. In this latter case $u$ is the barycenter of a
simplex $U^\prime$ that contains the vertex $w$, and hence again $u
\sim w$.  We conclude that $u \in N_{G_{k,X}^x}(w)$ and hence
$N_{G_{k,X}^x}(v) \subseteq N_{G_{k,X}^x}(w)$ as desired.

\begin{center}
\epsfig{file=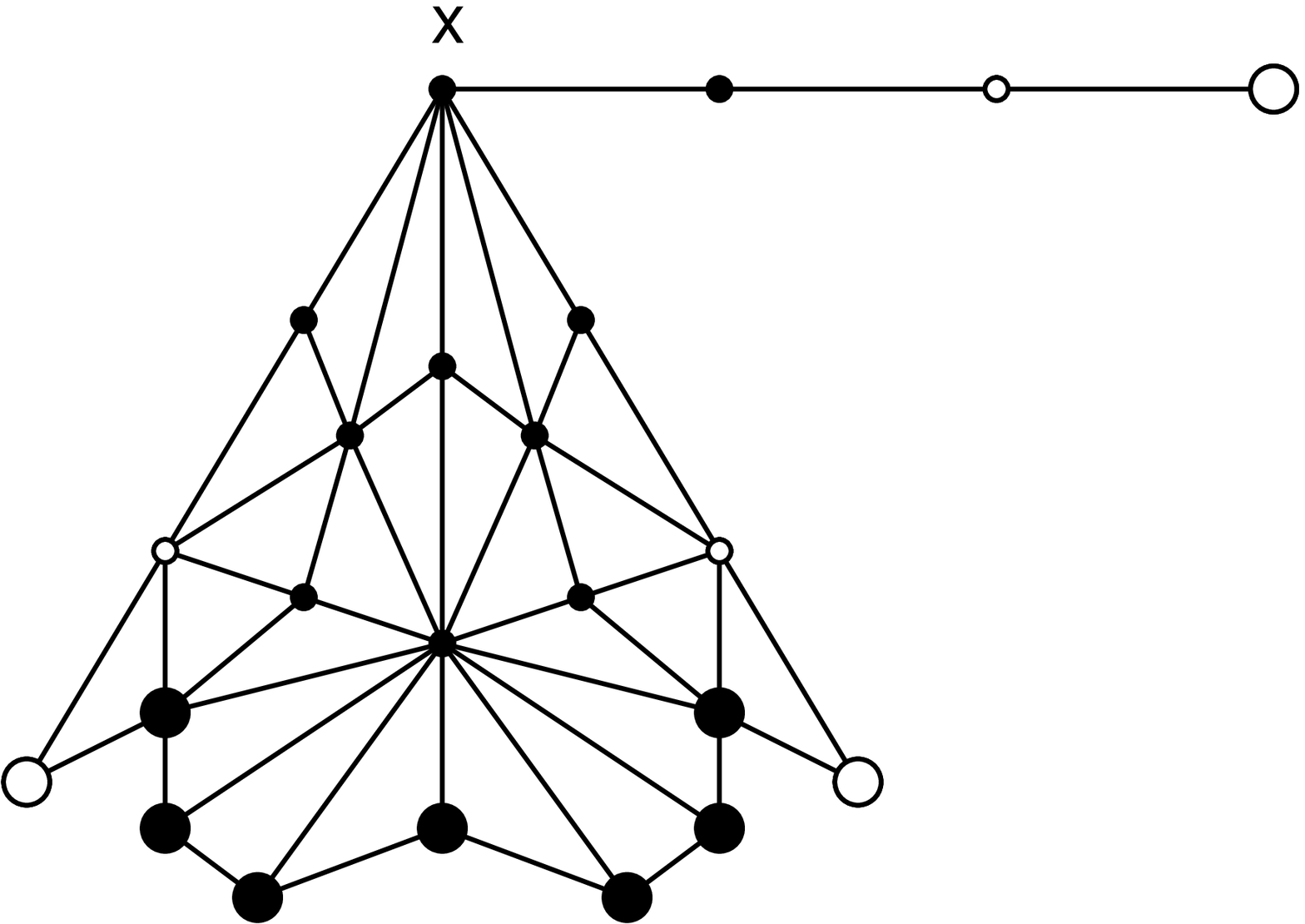, height=1.5 in, width = 2.4 in}
\epsfig{file=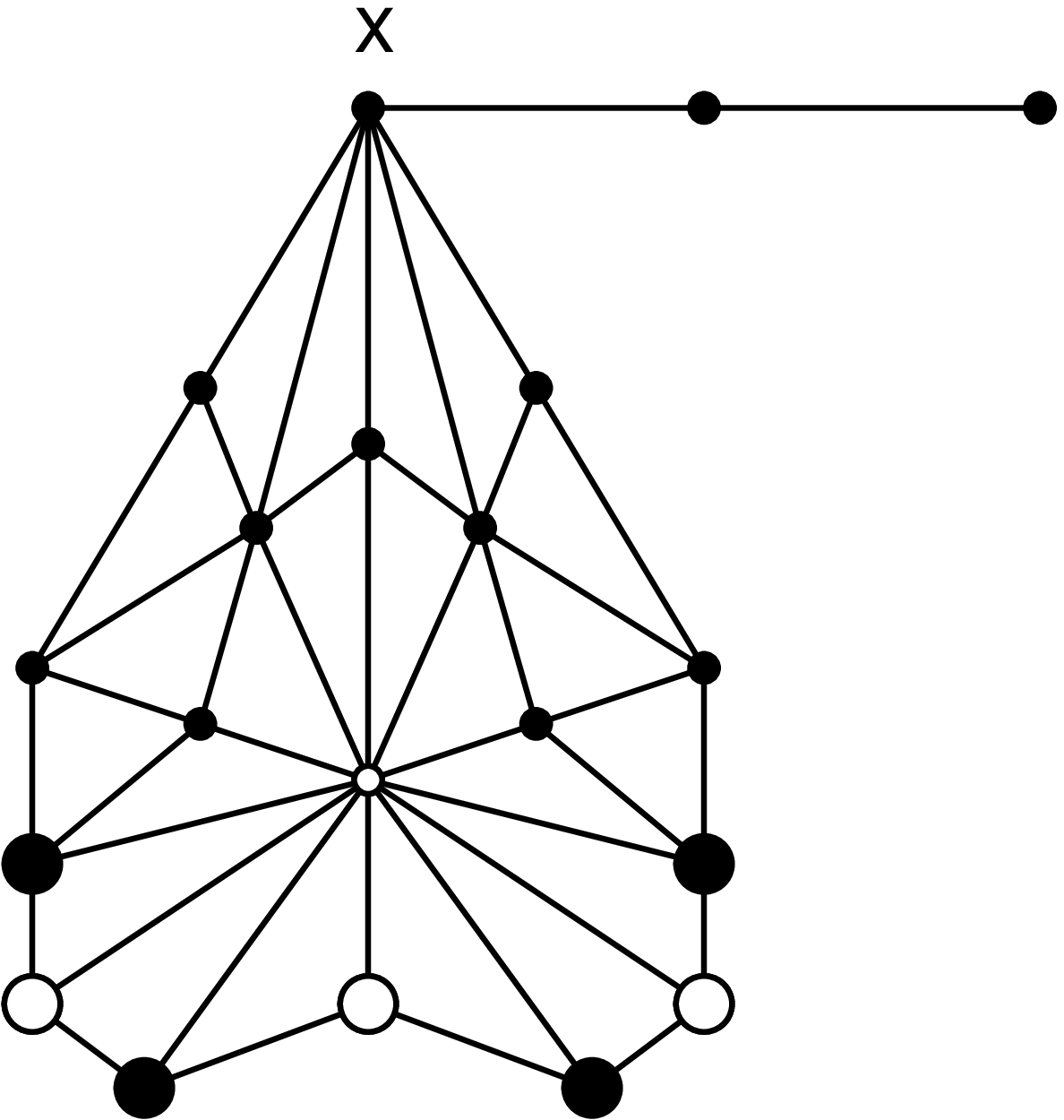, height=1.5 in, width = 1.65 in}

{Folding away vertices of type $(1,1)$ and of type $(2,1)$ in
$G^x_{2,X}$.}
\end{center}

\begin{center}
\epsfig{file=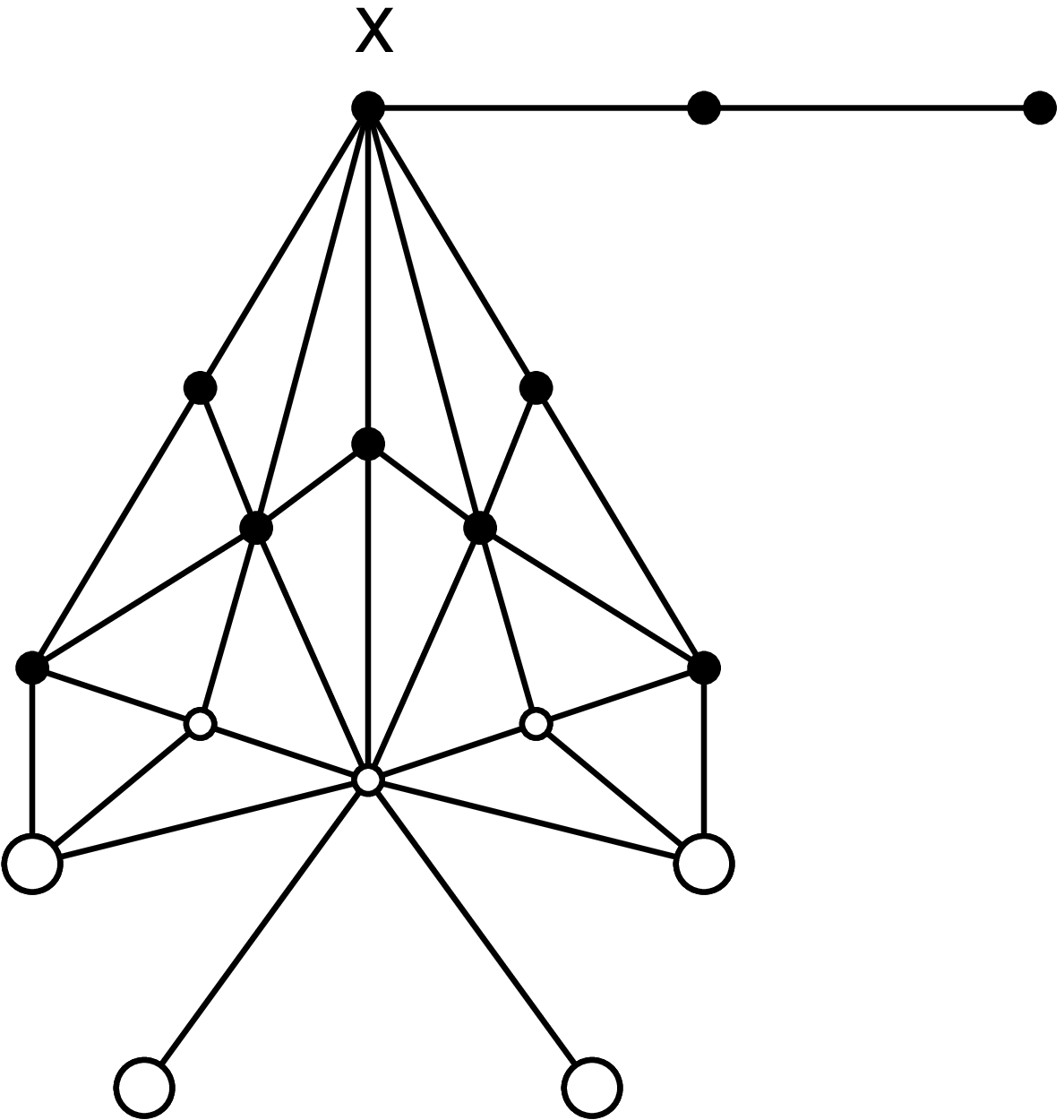, height=1.5 in, width = 1.65 in}
\epsfig{file=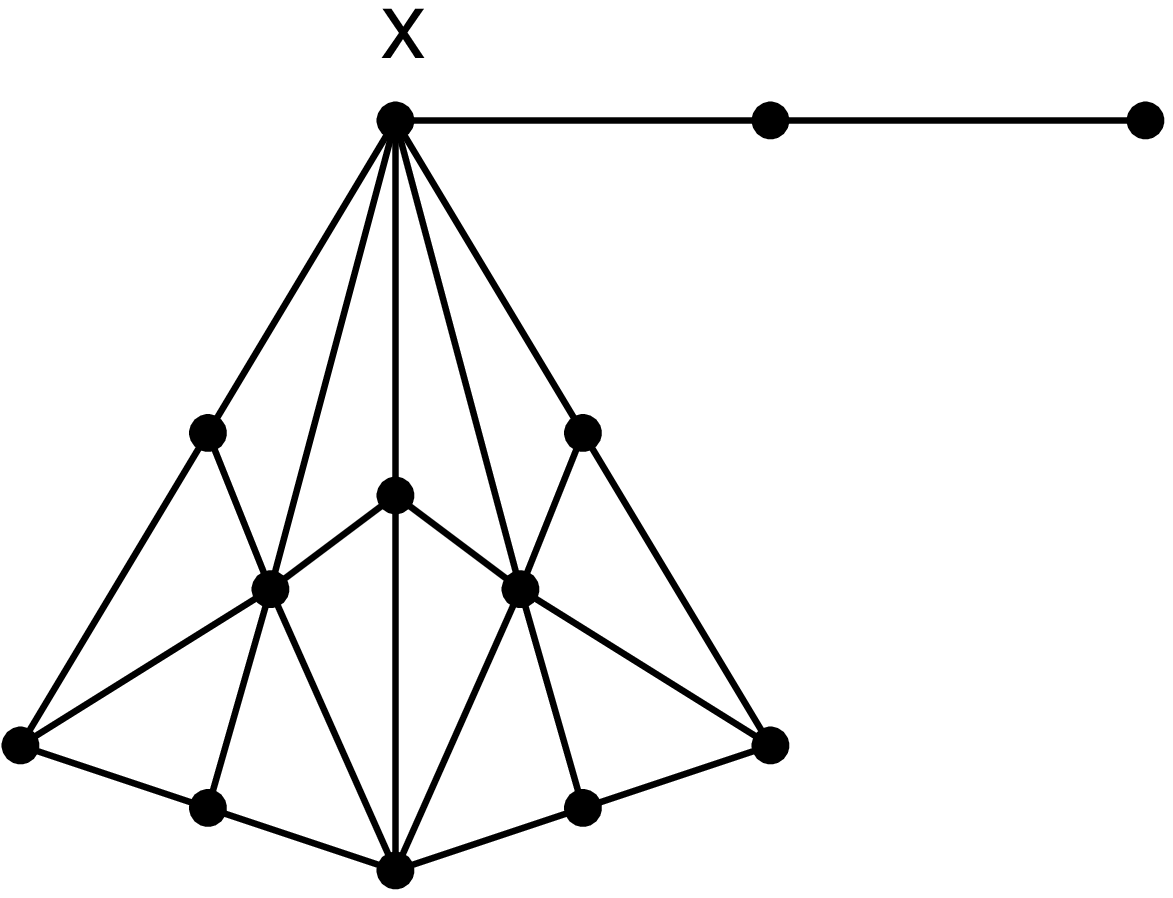, height=1.5 in, width = 1.65 in}

{Folding away vertices of type $(2,2)$ and the resulting
$(G^x_{1,X})^\prime$}.
\end{center}

This completes the induction on $j$ and hence we have now folded
away all vertices of $V_x$ that are of type $(i,?)$.  This in turn
completes the induction on $i$ and we conclude that all vertices in
$V_x$ can be folded away.  As we noted above, the resulting graph is
$\big(G_{k-1,X}^x \big)^\prime$, the barycentric subdivision of
$G_{k-1,X}^x$, which we conclude is dismantlable by induction on $k$
and by applying Lemma ~\ref{subdivision}.  The result follows.

\end{proof}

The final step in proving our theorem is to consider the
intersections of the subcomplexes $G_{k,X}^x$.

\begin{lemma} ~\label{intersection}
All nonempty intersections of the subcomplexes
$\big\{(G_{k,X}^x)^T\big\}_{x \in V(X)}$ are contractible.
\end{lemma}

\begin{proof}
We prove this in much the same way as we handled the contractibility
of the subcomplexes themselves.  In particular it is enough to show
that the subgraphs obtained as nonempty intersections of
$\{G_{k,X}^x\}_{x \in V(X)}$ are dismantlable.  A vertex of such a
graph is, by definition, within a distance of $2^k - 1$ of every
vertex $x \in V(X)$ in some index set $I \subseteq V(X)$.

Suppose $G_{k,X}^I$ is such a graph.  Again, we will show that
$G_{k,X}^I$ is dismantlable by induction on $k$.  If $k = 1$ then
the graph $G_{k,X}^I$ is a single looped vertex, the barycenter of
the face of $X$ defined by the index set $I$, which is of course
dismantlable.

For the case $k > 1$ we will, as above, fold away the vertices of
$G_{k,X}^I$ that are distance $2^k - 1$ from some vertex $x \in I$.
We will refer to these vertices as $V_I$, so that $V_I = \{v \in
G_{k,X}^I: d(v,x) = 2^k - 1$ for some $x \in I\}$.

Again, we fold away the vertices of $V_I$ in lexicographic order
according to their type $(i,j)$.  Since $V_I \subset V_x$ (for any
$x \in I)$, we can follow the same procedure as we described in the
proof of Lemma ~\ref{subcomplex}.  In particular, to fold away a
vertex $v \in V_I$ of type $(i,j)$, we choose a vertex $w \in
N_{V_x}(v)$ in the neighborhood of $v$ such that $w \in G_{k-1,X}^I$
and such that the type of $w$ is largest in the lexicographic order.

We just need to check that $w$ is within $2^k - 1$ of \textit{every}
vertex $x^\prime \in I$, so that indeed $w \in V_I$.  But this
follows from the choice of $w$: since $v$ is in the interior of the
face of $X$ determined by the vertices $I$, any neighbor $w^\prime$
of $v$ that lies outside of $V_I$ will be of type $(i^\prime,j)$,
where $i^\prime < i$.  But $v$ has neighbors in $G_{k-1,X}^I$ that
are of type $(i,j^\prime)$, so that the choice of $w$ will indeed
lie in $V_I$.

Hence the double induction follows through in this case, and we are
left with a graph $G_{k-1,X}^{I \prime}$, the barycentric
subdivision of the graph $G_{k-1,X}^I$ (informally speaking).  Once
again we employ Lemma ~\ref{subdivision} and by induction we get
that this graph is also dismantlable.

\end{proof}

\section{Further questions}
Having constructed our graph $G_{k,X}$ as the 1-skeleton of the
$k^{th}$ iterated subdivision of $X$, a natural question to ask is
if this choice of $k$ is best possible.  We have a feeling that it
is not, and in fact, for the case $\diam(T) = 1$ (so that $T$ is a
complete graph with possibly some loops) we conjecture that $k = 1$
will do the job.

\begin{conj}
If $X$ is a finite simplicial complex, and $T$ is a finite connected
graph with $\diam(T) = 1$, then there is a homotopy equivalence

\begin{center}
$\Hom(K_2, G_{1,X}) \simeq X$.
\end{center}

\end{conj}

Another thing to consider would be simplicial complexes with a
specified group action.

\begin{question}
Given a graph $T$ with automorphism group $\Gamma = Aut(T)$, and a
$\Gamma$-simplicial complex $X$, can one find a graph $G$ such that
$\Hom(T,G)$ is $\Gamma$-homotopy equivalent to $X$?
\end{question}

\section{Acknowledgements}
The author wishes to thank Carsten Schultz for fruitful discussions,
and especially his advisor, Eric Babson, who suggested the
construction of the graph $G_{k,X}$.

\bibliographystyle{halpha.bst}
\bibliography{litgraph}

\end{document}